\documentclass[12pt,a4paper]{article}

\usepackage{pdflscape}
\usepackage{amsmath}
\usepackage{amssymb}
\usepackage{latexsym}
\usepackage{srcltx}
\usepackage{graphics}
\usepackage{color}
\usepackage{epsfig}
\usepackage{color}

\newcommand{\co}{{\mathbb C}}

\newcommand{\re}{{\mathbb R}}
\newcommand{\n}{{\mathbb N}}

\newcommand{\cA}{{\mathcal{A}}}

\newcommand{\cE}{{\mathcal{E}}}

\newcommand{\cV}{{\mathcal{V}}}
\newcommand{\cL}{{\mathcal{L}}}

\newcommand{\cB}{{\mathcal{B}}}

\newcommand{\cS}{{\mathcal{S}}}
\newcommand{\cR}{{\mathcal{R}}}

\newcommand{\cP}{{\mathcal{P}}}
\newcommand{\cM}{{\mathcal{M}}}

\newcommand{\cC}{{\mathcal{C}}}

\newcommand{\bx}{{\boldsymbol{x}}}
\newcommand{\by}{{\boldsymbol{y}}}
\newcommand{\bz}{{\boldsymbol{z}}}

\newcommand{\bp}{{\boldsymbol{p}}}

\newcommand{\ba}{{\boldsymbol{a}}}

\newcommand{\bv}{{\boldsymbol{v}}}

\newcommand{\bm}{{\boldsymbol{m}}}

\newcommand{\bh}{{\boldsymbol{h}}}

\newtheorem{theorem}{Theorem}
\newtheorem{prop}{Proposition}
\newtheorem{lemma}{Lemma}
\newtheorem{cor}{Corollary}
\newtheorem{remark}{Remark}
\newtheorem{ex}{Example}
\newtheorem{defi}{Definition}

\date{}

\author{Vladimir Yu. Protasov 
\thanks{DISIM, University of L'Aquila, Italy; {e-mail: \tt\small
vladimir.protasov@univaq.it}} , 
Rinat Kamalov
\thanks{The Instutute of Control Sciences,  Russia {e-mail: \tt\small
rinat020398god@yandex.ru}} 
 }

\title{Stability of linear systems \\
with bounded switching intervals  
}

\begin{document}
\maketitle

\begin{abstract}

We address the stability problem for linear switching systems 
with mode-dependent restrictions on the  switching intervals. 
Their  lengths can be bounded as from below
(the guaranteed dwell-time) as from above. The upper bounds make this 
problem quite different from the classical case: a stable system can 
consist of unstable matrices, it may not possess Lyapunov functions, etc. 
We  introduce the concept of Lyapunov multifunction  with discrete monotonicity, 
which gives upper bounds for the Lyapunov exponent. Its existence as well as the existence 
of invariant norms are proved. Tight lower bounds are
 obtained in terms of a modified Berger-Wang formula over periodizable  
switching laws. Based on those results we develop a method of computation 
of the Lyapunov exponent with an arbitrary precision and analyse its efficiency in numerical results. The case when some of upper bounds can be cancelled is analysed.

\bigskip

\smallskip

\begin{flushright}

\end{flushright}

\end{abstract}
\bigskip

\vspace{1cm}

\begin{center}

\large{\textbf{1. Introduction}}	
\end{center}
\bigskip

Stability of continuous 
linear switching systems 
with  guaranteed mode-dependent dwell time, i.e., 
when the time interval between consecutive switches is bounded below,   has been studied in great detail 
due to applications  in industry such as robotics~\cite{K2010}, 
multilevel power converters~\cite{Basso}, etc. Theoretical study of such systems 
draw much attention in recent 
years~\cite{CGPS21, HespanhaMorse, Liberzon, LiberzonMorse, Morse}, see also~\cite{BS, Chesi1, chesi0, GC} and references therein for the numerical issues.  
In this work we address the stability problem when the time interval 
between switches can be bounded not only below but also above. 
Those are the systems that cannot stay in some modes too long and must be switched 
when the time interval reaches some fixed mode-dependent bound. This issue is 
natural both from theoretical and practical viewpoints, when a system 
can stagnate, overheat, etc., after being in one regime for a long time. 
Surprisingly enough, the approaches elaborated in the aforementioned works 
can hardly be applied to this problem, 
in spite of an apparent similarity. We begin with explanations of this phenomenon in Section~2.  
In fact, some properties of systems with upper bounds for the 
switching intervals are fundamentally different from both the classical systems 
(without the switching time restrictions) and the systems with guaranteed 
dwell time. First of all, those systems being stable may not have Lyapunov functions.  
They may not have Lyapunov multifunctions either, 
which are often associated to the guaranteed dwell time 
stability~\cite{CGPS21}. Second,  the stability of the system may not imply the 
stability of each matrix separately. Several unstable regimes can 
form a stable system,  although it may seem strange from the first site. 
More generally, a finite time switching law not always provides a lower bound 
for the Lyapunov exponent, as it is the case for usual systems. 
Thus, with the upper time restrictions we do not have the main source
of bounds for Lyapunov exponents. That is why  
systems with bounded switching time require new concepts of the very basic notions 
and actually developing a special theory. 

We deal with the linear switching systems of the form
\begin{equation}\label{eq.sys}  
\left\{
\begin{array}{l}
\dot \bx(t) \ = \ A(t) \bx(t), \qquad t \ge 0; \\
\bx(0) = \bx_0, \quad A(t) \, \in \, \cA\, . 
\end{array}
\right. 
\end{equation}
 Here~$\bx(\cdot ): \, \re_+ \to \re^d$ is the trajectory; the {\em switching law}  $A(\cdot )$
 is a measurable function with valuer on a compact 
 {\em control set}   of real~$d\times d$ matrices. We deal with the case of finite control set 
 $\cA \, = \, \{A_1, \ldots , A_n\}$. Each matrix~$A_j$ corresponds to its {\em mode}
 or {\em regime} of ODE~(\ref{eq.sys}). The time between two consecutive 
 switches is a {\em switching interval}. Each mode $A_j$ has restrictions on the length of its 
 switching interval between~$m_i$ (the lower bound) and~$M_j$ (the upper bound). 
 Thus, we obtain the {\em system with restrictions on the switching intervals}, or in short, 
 the {\em restricted system} 
 $S = \{\cA, \bm, \cM\}$, where $\bm = \{m_1, \ldots , m_n\}, \, \cM = \{M_1, \ldots , M_n\}$. 
 The system can stay in the regime~$A_j$ without switches  for the time interval of length on 
 $[m_j, M_j]$. 
 
 Note that if we impose the upper  bounds, then we have to restrict the dwell time as well.  Otherwise one could split every switching interval 
 by an arbitrary small interval of another regime, which would make 
 the upper restrictions useless. Thus, $m_i > 0$ and $M_i \in (m_i, +\infty]$
 for all~$i$. In the case~$M_i = +\infty$ there is no upper restriction
 for~$A_j$. In particular, if all~$M_j$  are infinite, then 
 we obtain a system with guaranteed mode-dependent dwell time. 
 If the converse is not stated we always assume that~$M_j$
 are finite. In Section~8 we show that this assumption is actually made 
 without loss of generality.

The restricted system is called {\em asymptotically stable}
if all its trajectories tend to zero as~$t\to +\infty$. 
We do not consider other kinds of stability and usually drop the word~``asymptotically''. 
For the classical (unrestricted) systems, the necessary condition of stability is that 
each matrix~$A_j$ is Hurwitz or, more generally, every finite switching law
(i.e., switching law on the finite segment) generates a Hurwitz matrix. 
Neither of this conditions work for the restricted systems. Finally, 
the usual stability condition provided by Lyapunov functions is not applicable either.
A stable restricted system does not necessarily possess Lyapunov function, 
i.e., a positively homogeneous function that decreases along every trajectory. 

To analyse the stability of restricted systems we begin with Lyapunov functions. 
The modification of this concept is done in several stages. 
First, we replace one function~$f$ but a family of functions~$\{f_j\}_{j=1}^n$
associated to the matrices~$A_j$. Each function~$f_j$ measures the trajectory~$\bx(t)$
on the interval of the regime~$A_j$. So, the Lyapunov function~$f$
possesses a piecewise structure, switching between the components~$f_j$
simultaneously with the modes~$A_j$. Second, we replace 
the  decrease condition of~$f(t)$ along any trajectory by the decrease 
along the sequence of the switching points. This discrete 
monotonicity ensures the stability of a restricted system. Moreover, 
we show that the existence of such a Lyapunov multifunction 
is necessary and sufficient for the stability. Every stable restricted system 
does possess a Lyapunov multifunction with convex components~$f_j$
(a {\em multinorm}) and, moreover, every irreducible system 
has an invariant multinorm, which is analogous to the Barabanov norm. 
To obtain lower bounds for the Lyapunov exponents, we prove a version of 
the Berger-Wang formula for restricted systems. While in the classical case, 
the Lyapunov exponent is estimated from below by finite switching laws, 
for restricted systems we use only {\em admissible} finite laws, which 
begin and end with different matrices. 

Applying the modified concepts of Lyapunov functions and of 
the Berger-Wang formula we derive an algorithm 
for computation of the Lyapunov exponent with an arbitrary precision 
and for deciding the stability of restricted systems. 
This algorithm originates with the invariant polytope algorithm
for unrestricted systems~\cite{GP13, GP16, Mej20} but it 
constructs simultaneously several polytopes as unit balls of the norms~$f_j$. 
The numerical results are demonstrated in Section~7.  

For the sake of simplicity we mostly deal with the case when 
all regimes~$A_j$ have the same lover and upper bounds~$m$ and $M$
respectively and denote such systems as~$S = \{\cA, m, M\}$.
All the results are extended to the mode-dependent restrictions~$m_j, M_j$ in a straightforward manner.  
In Section~2 we establish some basic properties of restricted systems and give  examples.
Section~3 introduces the modified concept of Lyapunov function as a 
multifunction with discrete monotonicity. In Section 4 we establish the 
fundamental theorems on the existence of the Lyapunov multifunction and of the invariant (Barabanov) 
multinorms. The latter concept is also modified for restricted systems. 
Section~6 presents the algorithm for computing the Lyapunov exponent 
and for deciding stability; the numerical results are provided in Section~7. 
Section~8 deals with the case when the upper bound~$M_j$ is big. 
Under which condition it can be reduced or, conversely, 
increased to~$+\infty$, i.e. be omitted? We solve this problem by introducing a concept of ``cut tail point'' of a matrix.  
Finally in Section~9 we provide proofs of the main result.

\begin{center}
\large{\textbf{2. Basic properties of restricted systems}}
\end{center}
\bigskip


We consider a linear switching system~$S = \{\cA, m, M\}$
defined by a family of matrices~$\{A_1, \ldots , A_n\}$
and by the time segment $[m, M]$ for the length of allowed length of 
switching intervals. Every switching law~$A(\cdot)$ is supposed to be admissible, 
i.e., having the lengths of all switching intervals on the segment~$[m, M]$. 
We assume formally that all the switching  intervals 
 are open from the right, so all the switching points $t_k$ are associated to the left 
 ends, i.e., 
 $A(t_k) \, = \, A(t_k - \varepsilon)$. The point~$t_0=0$ is formally associated 
 to some ``zero mode''~$A_k \in \cA$, so it is  
 considered as a ``zero switching point'' from~$A_k \in \cA$ to the first mode~$A_i$.

 Let us remember that we are restricted to the case 
 of equal time segments for all the modes for the sake of simplicity; 
 all our results are easily extended to arbitrary  mode-dependent segments~$[m_i, M_i]$.

 For an arbitrary switching law~$A(\cdot)$, we 
 consider the solution~$\Pi(t)$ of the matrix ODE
 $\Pi'(t) \, = A(t)\Pi(t), \, \Pi(0) = I$, where $I$ is the $d\times d$ identity matrix. 
 For every~$t\ge 0$ we call the~$P(t)$ a {\em product} generated by the switching 
 law~$A(\cdot)$, in analogy with the case of equal switching intervals of length~$\tau$, 
 when $\Pi(k\tau)$ is a product of matrices from the set~$e^{\tau A_j}, \, j = 1, \ldots , n$.
 Every trajectory  has the form~$\bx(t)\, = \, 
 \Pi(t)\bx_0$. A restriction of a switching low to a finite segment~$[0, T]$
 will be called a finite switching law. 
 
 Since the family~$\cA$ is compact, it follows easily that 
the family of functions~$\Pi(\cdot)$ is equicontinuous in each finite segment~$[0,T]$, 
and hence, by the Arzela-Ascoli theorem, it is compact on this segment. 
 Therefore, 
 the set of all trajectories~$\bx(t)$ starting 
in a given compact set~$K \subset \re^d$ is compact in~$C[0,T]$. 
This is true for all switching systems. For the restricted systems we can say more: 
 the set of admissible  switching laws is also compact. To formulate the corresponding property we need some further notation.  The set of all switching laws~$A(\cdot)$ of a restricted system on~$[0, T]$ will be referred to as~$\cE(T)$.  
This is a subset of the~$L_1$ space of martix-valued functions on~$[0,T]$. 
Similarly we define the set~$\cE = \cE(+\infty)$ with the topology 
of convergence in each segment in~$\re_+$. 
 \begin{prop}\label{p.5}
For arbitrary restricted system, the set of 
admissible switching laws~$\cE$  is compact and 
the sets~$\cE(T)$ are compact for all~$T>0$. 
\end{prop}
{\tt Proof.} Every switching law on~$[0,T]$ can have at most 
$T/m$ switching points, from which it depends continuously. Hence, 
the set~$\cE(T)$ is compact as an image of a compact set 
under a continuous map. Thus,~$\cE(T)$ is compact for all~$T$, hence, 
the set $\cE$ is compact as well.

{\hfill $\Box$}
\medskip

Similarly to the usual (unrestricted) switching systems,   the maximal rate of growth of trajectories is measured by the  Lyapunov exponent: 
  $$
 \sigma(S) \quad  = \quad \inf\ \Bigl\{\quad \alpha \quad \Bigl|\quad  
\|\bx(t)\| \ \le \ C\, e^{\, \alpha t} \quad \mbox{for every trajectory}\ \bx(t)\, \Bigr\}\, .
$$ 
For a system generated by one matrix~$A$, we have $\sigma(A)\, = \, \max_{k=1}^{d} {\rm Re}\, \lambda_k$ is the {\em spectral abscissa}, where ${\rm Re}$ denotes the real part of a complex number and 
 $\lambda_1, \ldots , \lambda_d$ are the eigenvalues of~$A$
counting multiplicities. In particular, the one-regime system is stable  if
the matrix~$A$ is {\em stable} or Hurwitz,  i.e., ${\rm Re}\, \lambda_k < 0$
for all~$k = 1, \ldots , d$.

The following fact is well known for unrestricted systems, however,  
in the restricted case it needs a different proof, which is given in Appendix.
\begin{prop}\label{p.10}
The system is stable precisely when~$\sigma(S)< 0$. 
\end{prop}
The Lyapunov exponent is usually estimated from above by 
a finite switching law. In particular, $\sigma(\cA) \ge \sigma(A_j)$
for all~$j=1, \ldots , n$. This is true for unrestricted systems 
or for systems with guaranteed dwell time~\cite{CGPS21}. However, if the system is restricted, then this inequality may fail, as the following simple example demonstrates. 

\begin{ex}\label{ex.5}(Stable system of unstable matrices). 
{\em Consider the system in $\re^2$ with two diagonal matrices
$A_1 \, = \, {\rm diag}\, (1, -3), A_2 \, = \, {\rm diag}\, (-3, 1)$
and with the interval $[m, M] = [1, 2]$. Clearly, both of them are unstable
with $\sigma = 1$. Nevertheless, the system is stable and~$\sigma\{A_1, A_2\} = -\frac13$. 
To see this we take an arbitrary switching law and denote by  $s_j \in [1, 2]$ the length of $j$th time interval.  
After the first $2k$ intervals ($2k-1$ switches), 
we obtain a diagonal  matrix~$\Pi_{2k}$
with the diagonal elements 
$e^{\ \sum_{j=1}^{k}\, (s_{2j-1} - 3s_{2j})}  \, \le \, e^{\, k(2 - 3) } \, 
= \, e^{\, -k }$ and 
$e^{\ \sum_{j=1}^{k}\, (-3s_{2j-1} + s_{2j})} \, \le \,  e^{\, -k}$. 
Since the total time is at most $3k$, we get~$\sigma(S) \le \frac13$. 
On the other hand, the switching law with $s_{2j-1} = 1, s_{2j} = 2, \, j \, \in \, \n$, 
gives the exponent of growth precisely~$\, \frac13$. 

Actually, this phenomenon is not surprising: the stationary switching laws 
$A(t) \equiv A_1$ and $A(t) \equiv A_2$ are unstable (the trajectory is unbounded)
but they both are not admissible.
}
\end{ex}
\smallskip 

The standard  way to prove the stability is 
to present a Lyapunov function, i.e., a positively homogeneous function 
decreasing along every trajectory of the system. It is well-known that 
every stable linear switching system does possess a 
Lyapunov function, which, moreover, can be chosen 
 convex (the Lyapunov norm)~\cite{MP1}. However, this is not true 
 for restricted systems. 

\begin{ex}\label{ex.10}(A stable restricted system which does not possess Lyapunov function). 
{\em Consider the system on $\re^1$ with two matrices (numbers)
$A_1 = a, \, A_2 = b, \ a > b$,  and with the switching interval~$[m, M] = [1,2]$. 
The maximal growth is given by the periodic switching law with the period~$3$:  $A(t) = a, \, 
t \in [0,2);  \ A(t)\, = \, b, \   t \in [2,3)$. The Lyapunov exponent is 
$\sigma \, = \, \frac{1}{3}\, ( 2a + b)$. Hence, if $2a + b < 0$, then the system is stable. 
For example, it is stable for  $a=1, b=-3$. However, there is no Lyapunov function in this case. 
Indeed, every Lyapunov function on~$\re$ has the form $f(\bx) \, = \, c\, |\bx|$
with some~$c>0$. Hence, if $0 = t_0 < t_1 < \cdots $ are the switching points, 
then in all the intervals~$(t_{2k}, t_{2k+1}), \, k\in \n$, we have $A(t) = a$, and hence 
$f\bigl(\bx(t)\bigr)\, = \, c\, |\bx(t)|\, = \, c\, |\, e^{\, a(t-t_{2k})} \bx(t_{2k})|$ increases on those intervals.  

Thus, the system is stable but $f$ somewhere increases, hence,  this is not a Lyapunov function.  
Let us remark that in this case even the sequence~$f\bigl(\bx(t_k)\bigr)$ at the switching points~$t_k$  is not monotone.
}
\end{ex}
\bigskip 

\begin{center}
\large{\textbf{3. The concept of Lyapunov functions for restricted systems}}
\end{center}
\bigskip 

Example~\ref{ex.10} above shows that a stable restricted system may not possess 
Lyapunov functions. In this case, how to prove stability? 
One needs to modify the the concept of Lyapunov function for restricted systems. 
We do this in three steps. First, we pass to multifunctions, which is a collection of functions: 
each  matrix~$A_j$ is equipped with its function~$f_j$. Second, the monotonicity 
requirement  for Lyapunov functions is relaxed to monotonicity over the 
set of switching points. Third, the Lyapunov multifunction is not defined 
in the minimal time intervals~$(t_k, t_k+m]$
 \medskip

 \begin{center}
\textbf{3.1 Three steps of modification}
\end{center}
\bigskip 

\textbf{Step 1. Lyapunov multifunctions}. The main idea is to 
consider the  action of each operator~$A_j$ separately of others.  
We use $n$ functions~$f_j$, a priory different, 
each of them corresponds to its operator~$A_j$ and measures the trajectory
$\bx(t)$ only on intervals of the regime~$A_j$. 
Thus, {\em the functions~$f_j$ are switched together with the operators~$A_j$}.
Now we compose the multifuntion~$f(\bx) = \bigl(f_1(\bx), \ldots , f_n(\bx)\bigr)$ with the following property:
$f(\bx(t))=f_j(\bx(t))$ when $A(t) = A_j, \ j=1,\ldots , n$. At the 
switching points the function~$f$ is supposed to be continuous from the 
left, i.e., if at the point $t_k$ the mode switches from~$A_i$ to $A_j$, 
then $f(\bx(t))=f_i(\bx(t))$. The point~$t=0$ can be associated to an arbitrary `` zero mode'' from~$\cA$.  

\smallskip 

 \textbf{Step 2. Asymptotic behaviour over  a discrete set}. 
Multifunctions do not solve the problem either: for a stable restricted system, there may  not exist functions~$f_j(\bx)$ that decrease 
on the intervals of the corresponding regimes~$A_j$. In Example~\ref{ex.10}, 
whatever positively homogeneous 
function~$f_j$ to take, it increases on each interval of the regime~$A_1 = a$. 

Thus, for restricted systems, one cannot provide the main 
property of Lyapunov functions to decrease along all trajectories. 
A way out of this situation is to narrow down the base of convergence. 
Namely, we relax the requirement for the  Lyapunov function to decrease
over the set of switching points only. 

For the restricted systems, when the lengths of switching intervals are bounded above,  
the asymptotics of an arbitrary trajectory $\bx(t)$ as $t\to \infty$
is equivalent to its asymptotics over the sequence of switching points. 
Hence, if a positive homogeneous function~$f$ 
decreases over switching points of every trajectory, then the system 
is stable, see Proposition~\ref{p.25}.

\smallskip 

\textbf{Step 3. The trajectory on the dwell time intervals}. 
We will 
pay no attention to the behaviour of the trajectory on the 
obligatory dwell time intervals~$(t_k, t_k+m]$, where we anyway cannot change the control. 
Due to the fixed lengths of those intervals, the asymptitic behaviour of the trajectory on them is the same as along the switching points. 
\smallskip

Let us now summarise the three main steps of the modification: 
\smallskip 

\noindent \textbf{1)} replace the concept of Lyapunov function~$f: \re^d \to \re_+$ by the Lyapunov multifunction\linebreak $f=(f_1, \ldots , f_n):\, \re^d \to \re_+^n$
and switch the functions~$f_j$ accordingly to the control regime simultaneously with the operators~$A_j$; 

\smallskip 

\noindent \textbf{2)} relax the condition for decreasing~$f(\bx(t))$ along the axis~$\re_+$ to
the decreasing on the set of switching points, i.e., 
$f(\bx(t_{k+1}))\,  <                                                                                                                                                                                                                                                                                                                                                                                                                                                                                                                                                                                                                                                                                                                                                                                                                                                                                                                                                                                                                                                                                                                                                                                                                                                                                                                                                                                                                                                                                                                                                                                                                                                                                                                                                                                                                                                                                                                                                                                                                                                                                                                                                                                                                                                                                                                                                                                                                                                                                                                                                                                                                                                                                                                                                                                                                                                                                                                                                                                                                                                                                                                                                                                                                                                                                                                                                                                                                                                                                                                                                                                                                                                                                                                                                                                                                                                                                                                                                                                                                                                                                                                                                                                                                                                                                                                                                                                                                                                                                                                                                                                                                                                                                                                                                                                                                                                                                                                                                                                                                                                                                                                                                                                                                                                                                                                                                                                                                                                                                                                                                                                                                                                                                                                                                                                                                                                                                                                                                                                                                                                                                                                                                                                                                                                                                                                                                                                                                                                                                                                                                                                                                                                                                                                                                                                                                                                                                                                                                                                                                                                                                                                                                                                                                                                                                                                                                                                                                                                                                                                                                                                                                                                                                                                                                                                                                                                                                                                                                                                                                                                                                                                                                                                                                                                                                                                                                                                                                                                                                                                                                                                                                                                                                                                                                                                                                                                                                                                                                                                                                                                                                                                                                                                                                                                                                                                                                                                                                                                                                                                                                                                                                                                                                                                                                                                                                                                                                                                                                                                                                                                                                                                                                                                                                 \, f(\bx(t_{k})), \ k \in \n$; 
\smallskip 

\noindent \textbf{3)}  do not consider the values of $f$ on the dwell-time intervals~$(t_k, t_k+m]$. 
 \smallskip

We are going to see that it is sufficient to consider only multinorms, 
when all~$f_j$ are convex functions. 
\begin{defi}\label{d.20}
The multinorm $f$ is called Lyapunov for a system~$\{\cA, m, M\}$ if 
the sequence $f\bigl(\bx (t_k) \bigr), \, k \ge 0,$ decreases on every trajectory~$\bx(t)$. 
\end{defi}
This means that if~$t_k$ is a switching point 
from the mode~$A_i$ to~$A_j$, then $f_{i}(t_{k}) \, > \, f_{j}(t_{k+1})$. 
The point~$t_0=0$ is considered as a ``zero switching point''
from  an arbitrary ``zero mode''~$A_k \in \cA$ to the mode starting mode~$A_i$. 
\begin{prop}\label{p.25}
If a restricted system~$S = \{\cA, m, M\}$ possesses the Lyapunov multinorm, then 
$S$ is stable. 
\end{prop}
{\tt Proof.} Let $f$ be a Lyapunov multinorm; 
then for each indices $i\ne j$, for every~$\bx \in \re^d$ such that $f_i(\bx)=1$
and for every~$\tau \in [m, M]$, we have $f_j(e^{\, \tau A_j}\bx) < 1$.  
Indeed, consider the trajectory~$\bx(\cdot)$ starting at
$\bx(0) \, = \, \bx$, with the formal zero mode~$A_i$ followed by the first 
mode~$A_j$ and $t_1=\tau$ is the first switching point. 
We have~$f_j\bigl(e^{\, \tau A_j}\bx\bigr) \, = \, 
f_j\bigl(\bx(t_1)\bigr)\, < \, f_i\bigl(\bx(t_0)\bigr) \, = \, 1$. 
Now by the compactness argument it follows that 
the maximum of $f_j\bigl(e^{\, \tau A_j}\bx\bigr)$ over all 
$\tau \in [m, M], \, i, j \in \{1, \ldots , n\}, \, i\ne j$
and over~$\bx$ such that~$f_i(\bx) = 1$ is less than one. 
Denote this maximum by~$q< 1$. Then for every trajectory~$\bx(\cdot)$
with~$f(\bx(0)) = 1$, 
we have $f(\bx(t_k))\, \le \, q^k$ and therefore, for every~$t \in [t_k, t_{k+1}]$, 
we have $f(\bx(t))\, \le \, C\, q^k$, where~$C$ is some constant. 
Therefore, $f(\bx(t))\to 0$ as~$t\to +\infty$.

{\hfill $\Box$}
\medskip 

\begin{remark}\label{r.18}
{\em Thus, if a positive homogeneous multifunction~$f$ 
decreases over switching points of every trajectory, then the system 
is stable.  For unrestricted systems this may not be true, and therefore, 
our modification of the Lyapunov function property is impossible. 

For the Lyapunov exponent of a restricted system, we have 
$$
\sigma(S)\ = \  
\inf\, \Bigl\{ \alpha \in \re \ \Bigl| \ \|\bx(t_k)\|\,  \le \, C\, e^{\, \alpha t} \, , \ k \in \n\Bigr\},
$$ 
where $\{t_k\}_{k\in \n}$ are the switching points,  and  the infimum is computed over all trajectories of the system.  
  Thus, for the  restricted systems, the asymptotics of 
any trajectory along the time axis is equivalent to that 
along the sequence of switching points. }
\end{remark}

One may wonder if the multinorm concept is really needed 
to have the Lyapunov property from Definition~\ref{d.20}? 
Can the components~$f_j$ be always chosen equal or at least proportional to each other? 
The answer is negative as the following example demonstrates.

\begin{ex}\label{ex.20}(The stable system which does not 
have a Lyapunov function with proportional components). 
{\em Consider the system on $\re^2$ with two matrices
$A_1, A_2$ defined as follows.  Take a diagonal matrix~$B \, = \, {\rm diag}\, 
\bigl( \frac{1}{2}, \frac{1}{3} \bigr)$ and a rotation~$R$ on an 
irrational ${\rm mod}\, \pi$ angle. Then $A_1, A_2$ are defined from the equalities  
$e^{A_1} \, = \, B ,  e^{A_1} \, = \, RB^{-1}$. 
For $m=M=1$ (which is formally not allowed but we make an exception for a while), 
each trajectory $\bx(t)$ is bounded. Furthermore,  $\bx(t_{2k})\, = \, \bx(2k)\, $ is an image of 
$\bx(t_0)$ under the rotation by the angle~$k\alpha$. 
Therefore, the points~$\bx(t_{2k})$ are everywhere dense on the circle. 
Hence, the Euclidean norm is the only norm which is non-decreasing on the 
sequence~$\{t_{2k}\}_{k\in \n}$. Assume $\|\bx(t_0)\| = 1$ and so, 
$\|\bx(t_{2k})\| = 1$ for all~$k$.  However, for the points~$\bx(t_{2k+1}) \, = \, 
B^{-1}R^k\bx(t_{0})$, we have $\|\bx(t_{2k+1})\| > 1$. 
Now take a close system  which is stable and 
$m, M$ are close but not equal. 
Thus, if this system has a Lyapunov multifunction $f=(f_1, f_2)$
with proportional $f_1, f_2$, then they both close to be proportional to the Euclidean norm, 
and hence $\|\bx(t_{2k+1})\| \, > \, \|\bx(t_{2k})\|$ for not very large~$k$.
Hence~$f$ cannot be a Lyapunov multinorm. We see that in this example there are 
no Lyapunov multinorms with proportional components.  
}
\end{ex}

Thus, to characterize stable systems one needs multinorms 
with a priori different (even non proportional) components.

To deal with Lyapunov multinorms it is convenient to 
to use  the graph interpretation. 
In the next subsection we recall the concept 
of dynamical systems on graphs and show how to put an arbitrary restricted switching system 
to the graph.  
\bigskip

\begin{center}
\textbf{3.2. The graph structure}
\end{center}
\bigskip 

Lyapunov multinorms  for restricted systems can naturally be formulated in terms of graphs. 
The switching systems of graphs have been  actively studied in the recent literature~\cite{CGP,K,PhJ1,PhJ3,SFS} and such models are special occurrences of hybrid 
automata~\cite{VdSS}. Some generalizations
to discrete-continuous hybrid  systems were presented in~\cite{CGPS21}, where 
they provided tools for the study systems with the lower restrictions on switching time. 
To the upper restrictions, that approach is hardly applicable, however, 
as we shall see below, the graph idea can also be put to good use. 
\medskip 

 \textbf{The graph multinorm}.
 We consider a directed graph~$G$ with vertices~$g_1, \ldots , g_n$, each vertex~$g_j$ is 
 associated to a linear space~$L_j$ isomorphic to~$\re^d$. Let 
 $\cB = \{B_1, \ldots , B_n\}$ be a family of linear operators on~$\re^d$. 
 All edges of~$G$ are coloured 
 in $n$ colours corresponding to these operators. 
 Each pair of vertices~$g_i, g_j, \, i\ne j$, are connected by two edges: 
 the edge of colour~$B_j$ is directed from $g_i$ to $g_j$ and 
 the opposite edge is of colour~$B_i$. Thus, each vertex~$g_j$ has precisely 
 $n$ incoming edges, all of the colour~$B_j$. The operator~$B_j$ 
 associated to the edge~$g_i g_j$ acts from~$L_i$ to $L_j$
 (so~$B_j: \re^d \to \re^d$ is considered as an operator from~$L_i$ to $L_j$). 
Now we introduce a multinorm on this graph. 
\begin{defi}\label{d.10}
If every  space $L_i$ on the graph $G$ is equipped with a norm $\|\cdot \|_i$, then
the collection of norms $\, \bigl\{ \, \|\cdot \|_i, \ \bigl| \ 
i = 1, \ldots , n\, \bigr\}$, is called  a multinorm.
The norm of an operator $B: L_i \to L_j$ is defined as
$\, \|B\| \, = \sup\limits_{\bx \in L_i, \|\bx\|_i = 1}\|B\bx\|_j$.
\end{defi}
We denote that multinorm by $\|\cdot\|$. Thus, for every~$\bx \in L_i$, we 
have  $\|\bx\| = \|\bx\|_i$ and we drop the index of the norm if it is clear
to which space $L_i$ the point $\bx$ belongs to.
\medskip 

\textbf{From  restricted systems to dynamical systems on graph}. 
 We consider an arbitrary restricted system~$\{\cA, m, M\}$ in~$\re^d$ and turn it to a system on a graph~$G = \{g_i\}_{i=1}^n$ with $L_j = \re^d, j= 1, \ldots , n$. 
We assume that $A_j$ acts in the space~$L_j$ and 
denote~$B_j = e^{\,mA_j}, \, j = 1,\ldots , n$.     
\smallskip 

Define
the following linear switching system on the graph~$G$. 
We start with a point~$\bx(t_0) \in L_i$ (for some~$i$), 
then jump to some~$L_j$ to the point~$\bx(t_0+m)\, = \, B_j\bx(t_0)$, then 
go in~$L_j$ along the curve~$\bx(t) \, = \, e^{(t-t_0)A_j}\bx(t_0), \ 
t \in [t_0+m, t_1]$ with some~$t_1 \le t_0+M$. In other words, 
we go along the solution of the ODE $\ \dot \bx(t) \, = \, A_j\bx(t), \, t \in [t_0+m, t_1]$
from the point~$\bx(t_0+m)\, = \, B_j\bx(t_0)$ during the time at most~$M-m$
until the switching point~$\bx_1$. Then we apply some~$B_s$
and jump to the space~$L_s$ landing at the point~$\bx(t_1+m)\, = \, B_k\, \bx(t_1)$, etc. 
We obtain a trajectory~$\bx(t)$ with switching points~$\{t_k\}_{k\ge 0}$. 

In our notation, for each~$k\ge 0$, the point~$\bx(t_k)$ is in the space~$L_{j_k}$. From this point the 
trajectory goes to the space~$L_{j_{k+1}}$ by means of the 
operator~$B_{j_{k+1}}$. Thus, the trajectory switches from~$\bx(t_k)$
 to  the point~$\bx(t_k +m) \, = \, 
B_{j_{k+1}}\bx(t_k)$ in the space~$L_{j_{k+1}}$. Then the trajectory 
stays inside this space for the time~$t_{k+1} - t_k - m$
and goes along the  solution of the ODE $\, \dot \bx(t) = A_{j_{k+1}}\bx(t), \, t \in [t_k+m, t_{k+1}]$. 

We see that  each jump (switch) between the spaces takes time precisely~$m$ and each 
ODE takes time $t_{k+1} - t_k -m$, which does not exceed~$M-m$. 
We control the switching points~$t_k$ under the restrictions
$m \le t_{k+1}- t_k \le M$ for all~$k \ge 0$,  and the choice of the indices~$j_k$. 
Each choice produces a control law  and a trajectory~$\bx(t)$. 
For $t \in (t_k, t_k+m)$ the point $\bx(t)$ is not considered, although 
at can be defined as~$e^{(t-t_k)A_{j_{k+1}}}\bx(t_k)$.   
\smallskip 

\bigskip 

\begin{center}
\large{ \textbf{4. Lyapunov and invariant multinorms}}
\end{center}

\bigskip 

The Lyapunov multinorm  of a restricted switching system (Definition~\ref{d.20}) has an obvious geometrical interpretation on the graph. Let $\cB_i$ be the unit ball of the norm~$f_i$ and $\cS_i$ be the corresponding unit sphere. 
\begin{prop}\label{p.20}
A multinorm~$f$ is Lyapunov if for each $i = 1,\ldots , n$, and $j\ne i$, 
and for every~$\bx \in \cS_i$, the curve 
$e^{t A_{j}}\, B_{j}\, \bx \, , \ t \in [0, M-m]$,  
lies strictly inside~$\cS_{j}$. 
\end{prop}
{\tt Proof.} If $f$ is a Lyapunov norm, then 
choosing arbitrary $\bx(t_1) \in \cS_i$ and $t_2 \in [m, M]$
and denoting~$t=t_2-m$, we 
get $\bx(t_2)\, = \, e^{\, t A_{j}}\, B_{j}\, \bx(t_1)$, and therefore,  
${f\bigl(\bx(t_2) \bigr) \, < \, f\bigl(\bx(t_1) \bigr)\, = \, 1}$. Hence, 
the point $\bx(t_2)$ is strictly inside~$\cS_{j}$. Conversely, if the condition of the proposition is satisfied, then, for an arbitrary trajectory~$\bx(\cdot)$
and for an arbitrary its switching point~$t_k$ from a mode~$A_i$ to~$A_j$, 
we have $\bx(t_{k+1})\, = \, e^{t A_{j}}\, B_j\, \bx(t_k)$ with 
$t = t_{k+1} - t_k -m$. Normalising~$\bx(t_k)$ to have norm one, 
we obtain the point~$\bx(t_{k+1})$ inside~$\cS_j$, hence~$\|\bx(t_{k+1})\| < 1$. 
Thus, $\|\bx(t_{k+1})\| < \|\bx(t_{k})\|$, which completes the proof.

{\hfill $\Box$}
\medskip 

\begin{cor}\label{c.10}
If $f$ is a Lyapunov multinorm, then for every trajectory~$\bx$ 
and for every its switching point~$t_k$, we have $f\bigl( \bx(t) \bigr) \, < \, f\bigl( \bx(t_k) \bigr)$ for all~$t > t_k$. 
\end{cor}
Thus, along  every trajectory~$\bx(\cdot)$,  the Lyapunov function ``decreases by intervals'': 
at each point of the interval $(t_{k}, t_{k+1})$ it is less than at each point of  the previous 
interval~$(t_{k-1}, t_{k})$.

The following theorem shows that the existence of a Lyapunov multinorm with discrete monotonicity completely 
characterizes stable systems. 
\begin{theorem}\label{th.20}
Every stable restricted system possesses a Lyapunov multinorm. 
\end{theorem}
The proof is  in Section~9. The next stage is the extension of the concept 
of invariant or Barabanov norms to restricted systems. Such norms are usually applied 
to obtain a very refined analysis of the fastest growth of trajectories. 
Now it is quite expected that such an extension should require 
significant modifications for restricted systems.

\begin{defi}\label{d.30}
An invariant  multinorm for a system~$S = \{\cA, m, M\}$
is a multinorm~$f$ with two properties: 
\smallskip 

 \noindent \textbf{a)} 
for every trajectory~$\bx(\cdot )$ with switching points~$\{t_k\}$, 
the sequence $e^{-\sigma\, t_k}\, f\bigl(\bx (t_k) \bigr)$ is non-increasing; 
 \smallskip 
 
 \noindent \textbf{b)} 
for every $\bx_0 \in \re^d$, there exists a trajectory~$\bar \bx(\cdot )$ starting at~$\bx_0$ such that 
$e^{-\sigma\, t_k} f\bigl(\bar \bx (t_k) \bigr)\, \equiv \, {\rm const}$.  
\end{defi}
\begin{remark}\label{r.40}
{\em If $\sigma(S) = 0$, then the definition of the invariant multinorm 
gets a simpler form: }

for every trajectory, the sequence
$f\bigl(\bx (t_k) \bigr)$ is non-increasing  and for every~$\bx_0\in \re^d$, 
there exists a trajectory~$\bar \bx(\cdot)$ starting at~$\bx_0$ such that~$f\bigl(\bar \bx (t_k) \bigr)\, \equiv \, {\rm const}$. 
\smallskip 

{\em Note that 
every system can be reduced to the case~$\sigma = 0$ by replacing each matrix~$A_i$
with~$A_i \, - \, \sigma I$. Hence, it is sufficient to formulate 
the notion of invariant multinorm only for the case~$\sigma = 0$.

}
\end{remark}

\begin{prop}\label{p.30}
Let $\{\cA, m, M\}$ be a system with~$\sigma = 0$. A multinorm~$f$ is invariant  if for each $i = 1,\ldots , n$
 and for every~$\bx \in \cS_i$, the curve 
$e^{t A_{i+1}}\, B_{i+1}\, \bx \, , \ t \in [0, M-m]$,  
lies in the ball~$\cB_{i+1}$ and at least one point of this curve 
is on its surface~$\cS_{i+1}$. 
\end{prop}
{\tt Proof.} If $f$ is invariant,  then 
choosing arbitrary $\bx(t_1) \in \cS_1$ and $t_2 \in [m, M]$
and denoting~$t=t_2-m$, we obtain 
$f\bigl(\bx(t_2) \bigr) \, \le \, f\bigl(\bx(t_1) \bigr)\, = \, 1$. 
Moreover, this inequality must become an equality for some 
$t_2 \in [m, M]$, in which case $\bx(t_2)$ is on the surface of~$\cS_{i+1}$. 

Conversely, arguing as above we conclude that 
for each choice of switching points~$\{t_k\}_{k\in \n}$, 
it holds that~$f\bigl(\bx (t_{k+1})\bigr) \le f\bigl(\bx (t_{k})\bigr), \ k \in \n $.
The trajectory~$\bar \bx$ is constructed by induction. 
The point~$\bar \bx(t_1) \in \cS_1$ is arbitrary. Let the trajectory~$\bar \bx$ be 
constructed on the segment~$[0, t_k]$. Applying our assumption to the point 
$\bx = \bar \bx(t_k)$ we get a point $e^{t A_{i+1}}\, B_{i+1}\, \bar \bx (t_k)$ on the unit sphere for some $t\in [0, M-n]$.   We set $\, t_{k+1} = t_k + t$ 
and obtain~$f\bigl( \bx(t_{k+1})\bigr) = 1$. This way we extend 
the trajectory~$\bar \bx(t)$ to the  segment~$[t_k, t_{k+1}]$,  which  
proves the inductive assumption. 

{\hfill $\Box$}
\medskip 

Every trajectory of the 
maximal growth takes values~$f(\bx(t_k))\, = \, e^{\, \sigma t_k}f(\bx(0))$ at all its 
switching points~$t_k, \, k \ge 0$, while for all other trajectories this sequence 
in non-increasing and for some~$k$ we have~$f(\bx(t_k))\, < \, e^{\, \sigma t_k}f(\bx(0))$  (for each trajectory, the values of~$f$ are computed over its own 
switching points). Thus, the invariant multinorm allows us to characterize all 
trajectories of maximal growth and to evaluate  the multinorm at all its switching points. 
A question arises if such a norm always exists? The answer is affirmative provided 
the family of matrices~$\cA$ is {\em irreducible}, i.e., those matrices do not share a 
common nontrivial invariant subspace. 
\begin{theorem}\label{th.30}
Every irreducible restricted system
possesses an invariant Lyapunov multinorm. 
\end{theorem}
The proof is in Section~9.   It is interesting that Theorem~\ref{th.30} 
has the same formulation as for the usual (unrestricted) systems, proved by Barabanov 
in~\cite{B} with one exception: for the usual system it holds 
provided the control set~$\cA$ is convex (and hence, infinite, unless one-point). 
This distinction is explained by the different definitions: the invariant norm according to Definition~\ref{d.30}  exists for all finite~$\cA$. 
 Another modification of invariant norm was done for hybrid systems in~\cite{CGPS21}. 
  Also note that 
for systems with at least two matrices, the irreducibility is a generic property. 
Moreover, it suffices to analyse the maximal growth of trajectories and find the
Lyapunov exponents for irreducible systems only (see Section~6).

\medskip

\bigskip 

\begin{center}
\large{\textbf{5. Periodic solutions. The 
Berger-Wang formula}}
\end{center}
\bigskip 

One more difference of the restricted systems from the classical ones is the role of periodic trajectories. 
In the classical case,  for the unrestricted systems, each finite switching law~$A(\cdot) \in \cE(T)$
gives a lower bound for the Lyapunov exponent:~$\sigma(T) \, = \, 
T^{-1}\sigma\bigl(\Pi(T)\bigr)$, where, recall, $\Pi(t)$ is 
the solution of matrix ODE $\Pi'(t) = A(t)\Pi(t), \, \Pi(0) = I$. 
 Indeed, the maximal growth over all trajectories 
does not exceed the growth of the periodic switching law with the period~$\Pi(T)$. 
For restricted systems this is not true. For instance,  the univariate system 
from Example~\ref{ex.10} has the Lyapunov exponent~$\sigma = \frac{2\ln a \, + \, \ln b}{3} $, while 
the stationary switching law $A(t) \equiv a$, which is, of course, also periodic, has 
the Lyapunov exponent $\ln a \, > \, \sigma$. The reason is that 
not all finite switching laws can be periodized  for a restricted switching system. 
The obtained periodic law must belong to~$\cE$, i.e., have the time intervals for each regime~$A_i$ are all between $m$ and~$M$.

\begin{center}
\textbf{5.1. Admissible periods}
\end{center}
\bigskip 

 Which finite switching laws are periodisable? 
These are the {\em admissible} laws according to the following definition:  
 \begin{defi}\label{d.45}
 A finite switching law is called admissible if it begins and ends with different matrices.
 The set of admissible switching laws is denoted as~$\cE_0$.  
 \end{defi}
\begin{prop}\label{p.40}
 A finite switching law belongs to~$\cE_0$
precisely when there exists~$i \in \{1, \ldots , n\}$ such that this 
law on the graph starts and finishes in~$L_i$. 
\end{prop}
{\tt Proof.} If~$A(\cdot) \in \cE_0$ starts with a matrix~$A_j$ and ends with~$A_i$, then 
$i\ne j$. Moreover, the time interval for~$A_i$ is at least $m$, hence any product 
$\Pi(t)$ starts with $e^{mA_j} = B_j$. Taking arbitrary~$\bx \in L_i$ we can apply~$B_j$
and obtain a point~$B_j\bx \in L_j$, then we continue along the switching law~$A(\cdot)$
and finally return to~$L_i$ since the law finishes with~$A_i$. 
Conversely, if a finite trajectory starts and finishes in~$L_i$, then 
the corresponding switching law~$A(\cdot)$ starts with some matrix different 
from~$A_i$ and finishes with~$A_i$. Hence, $A(\cdot) \in \cE_0$. 

{\hfill $\Box$}
\smallskip 

Proposition~\ref{p.40} yields  that $\cE_0$ is 
the set of all periodizable finite switching laws:
\begin{cor}\label{c.20}
A periodic switching law belongs to~$\cE$ if and only if its period is from~$\cE_0$. 
\end{cor}
Now we can show that every $A(\cdot)\in \cE_0$ provides a lower bound for the Lyapunov exponent.
Recall that for matrix~$\Pi$ we denote by $\sigma(\Pi)$ its spectral abscissa, 
i.e., the biggest real part of its eigenvalues. 
\begin{prop}\label{p.50}
Fro a restricted system~$S$, every product~$\Pi(T)$
generated by a switching law from  $\cE_0$ possesses the property  $\ T^{-1}\sigma(\Pi(T)) \, \le \, \sigma(S)$. 
\end{prop}
{\tt Proof.} Since~$A(\cdot) \in \cE_0$, there exists an admissible  periodic switching law 
with period~$\Pi(T)$.  The  rate of growth  of this switching law 
is~$T^{-1}\bigl(\sigma(\Pi(T)\bigr)$ and 
this does not exceed the maximal rate of growth, which is~$\sigma(S)$.   

{\hfill $\Box$}
\medskip

Thus,  finite switching laws from~$\cE_0$ provide lower bounds for~$\sigma(S)$. 
The question arises if they can be arbitrarily sharp? In other words, can the 
values~$\, T^{-1}\sigma(\Pi(T))$ for $\, A(\cdot )\in \cE_0$  approach arbitrarily close to~$\sigma(S)$?   
For the classical (unrestricted) systems the answer is affirmative
because all finite switching laws 
can be periodized. This implies the   
 classical Berger-Wang formula for the classical case:  
\begin{equation}\label{eq.BerW-stand}
\sigma(S) \quad = \quad \limsup_{T\to \infty}\ 
\frac{1}{T}\, \max\limits_{A(\cdot) \in \cE(T)}\, \sigma \bigl(\Pi(T)\bigr)\, .  
\end{equation}
An extension of this formula to restricted systems is not direct 
since in this case we have a narrower class of admissible 
finite laws~$E$ and, respectively, the less set of lower bounds. 
That is why there is no evidence that those bounds  reach~$\sigma(S)$. 
Nevertheless, this is true. The following  corollary of Theorem~\ref{th.30} 
not only generalizes the 
Berger-Wang formula to restricted switching systems but even slightly improve it. 

\begin{theorem}\label{th.40}
For every restricted switching system, we have 
\begin{equation}\label{eq.BerW}
T\, \sigma \, (\cA)\ - \ \max\limits_{A(\cdot) \in \cE_0(t), \, t \le T}\, 
\sigma (\Pi(t))\quad \to \quad 0\ , 
\qquad \mbox{as}\quad  T \to \infty\, . 
\end{equation}
\end{theorem}
The proof  is in Section~9. 
Dividing equality~(\ref{eq.BerW})  by~$T$, we obtain~(\ref{eq.BerW-stand}). 
\begin{cor}\label{c.40}
The equality~(\ref{eq.BerW-stand}) holds for restricted systems after replacing 
~$\cE$ by~$\cE_0$ . 
\end{cor}
For unrestricted systems, the assertion~(\ref{eq.BerW}) implies the  Berger-Wang formula
but not vice versa.  Hence, even in the classical case Theorem~\ref{th.40} strengthens the  Berger-Wang formula. 
\begin{defi}\label{d.70}
For given~$T>0$,~{\em the maximal $T$-period} of the restricted system 
is a switching law~$A(\cdot) \in \cE_0(t), \ t \in [m, T]$ for which the value 
$\, t^{-1}\, \sigma(\Pi(t))$ is maximal among all admissible switching laws in the 
the time interval no longer than~$T$. 
\end{defi}
The maximal~$T$-period is well-defined because the set of  products
$\Pi(t)$ generated by all switching laws~$A(\cdot) \in \cE_0(t), \,  \ t \in [m,T]$,    is compact.

The maximum in~Theorem~\ref{th.40} can be computed over the maximal $T$-periods. 
  It would be natural to expect some regular behaviour of the maximal  $T$-period 
  with respect to~$T$, for example, its continuity. Suppose that 
  the matrix~$A_j$ is dominant in $\cA$, i.e., 
  it has the maximal spectral abscissa~$\sigma(A)$ among all~$A \in \cA$. Is it true that 
 the maximal $(m+M)$-period is provided by 
 a law with the maximal contribution of~$A_j$, i.e., when 
the mode~$A_j$ has the time interval of length~$M$?

The answer is negative. Moreover, the maximal $T$-period can be discontinuous  in~$T$. 
Even its length can be discontinuous. 

\bigskip 

\begin{center}
\textbf{5.2. A restricted system with discontinuous $T$-period}
\end{center}
\bigskip

Consider the system~$\cA$ of two $2\times 2$ matrices with $m=\pi$ and with arbitrary $M> \pi$: 
$$
A_1\ = \ 
\left( 
\begin{array}{rr}
-1&0\\
0 & -a
\end{array}
\right)\ ; 
\qquad 
A_2\ = \ 
\left( 
\begin{array}{rr}
0&-1\\
1 & 0
\end{array}
\right)\ ,  
$$
where $a> 1$. We have   $e^{tA_1}\, = \, {\rm diag}\, (e^{-t}, e^{-2t})$
and $e^{tA_2}$ is a rotation by the angle~$t$. Therefore, for the 
unrestricted system with the matrices~$A_1, A_2$, the matrix~$A_2$ is dominant and 
the law~$A(t)\equiv A_2$ provides the maximal growth. Does $A_2$ preserve its 
dominancy for the constraint system? This means that 
the longer the time interval for~$A_2$ the faster growth we have. 
For a given~$T$, denote by~$\Pi_m(T)$ 
the product of the following switching law: ~$A(t) = A_1, \, t \in [0,m)$ and $A(t) = A_2, \, [t, T]$. The question is whether~$\Pi_m(M)$ is the $M$-period whatever 
upper time restriction~$M$ to take?  

Before giving the answer let us make some observations.  
First,  the only law~$A(\cdot)\in \cE_0$, 
for which the (Euclidean) norm $\|Pi(T)\|$ reaches the value~$e^{-\pi}$ is 
$\Pi_m$. Indeed, the rotation~$e^{tA_2}$ does not change the norm while the contraction~$e^{tA_2}$
reduce it. For $t=\pi k$ the matrix~$\Pi_m(t)$
has the spectral radius~$e^{-\pi}$ (the corresponding eigenvector is~$\bv = (1,0)$). 
Hence, the spectral radius of~$\Pi_m(t)$ reaches its 
maximal value~$e^{-\pi}$ precisely at the points~$t = \pi k$. 
For all other laws, the spectral radius of~$\Pi(t)$ is smaller than~$e^{-\pi}$
for all~$t$. 

\begin{prop}\label{p.70}
For every large $a>1$, there is a number~$n\in \n$ and positive numbers
$\tau_1 < \cdots < \tau_{2n}$ such that for each~$k = 1, \ldots , n$, 
the following hold:

1) $\ \pi k < t_{2k-1} < \pi (k+1)$; 

2) for the 
restricted system~$\cA$ with~$M=\tau_{2k}$, a the solution~$E_m(M)$ is a unique~$M$-period;   

3) for the restricted system~$\cA$ with~$M=\tau_{2k-1}$. 

\noindent Moreover, $n\to \infty$ as $p\to \infty$. 
\end{prop}
{\tt Proof.}
By the compactness, for each~$t$, 
the maximum of~$\rho(\Pi(t))$ over all laws~$A(\cdot)\in \cE_0$ 
is achieved  and is smaller than~$e^{-\pi}$ for all points~$t\ne \pi k,\, k \in \n$. 
Denote by~$\varphi(t)$ this maximum to the power~$1/t$ and 
by $\varphi_m(t)$ the spectral radius of~$\Pi_m(t)$
to the power~$1/t$. We have $\varphi(\pi k) = \varphi_m(\pi k)\, = \, e^{-1/k}$
and $\varphi(t) \ge \varphi_m(t)$ for all other~$t$. 
 Observe that for the spectral radius of~$\Pi_m\bigl(\frac{\pi}{2} + \pi k \bigr)$
 is equal to~$e^{-\frac{(a+1) \pi }{2}}$, therefore 
 $$
 \varphi_m\, \left(\frac{\pi}{2} + \pi k  \right)\ = \ 
e^{-\frac{(a+1)}{2k+1}} \ < \ e^{-\frac{1}{k}}\ = \ 
\varphi_m\, \left(\pi k  \right)\, . 
 $$
 Thus, for each segment~$\bigl[\pi k\, , \, \pi (k+1)\bigr]$, the value of~$\varphi_m$
 at the midpoint is less than at the ends.    Denote by~$h_{2k-1}$ the minimum of~$\varphi_m$ on  this  segment and by~$h_{2k}$ its maximum on the segment~$\bigl[\, \pi \bigr(k + \frac12\bigr)\, , \, \pi \, \bigl(k + \frac32\bigr)\bigr]$. Then~$h_{2k}$
 is the point of maximum of~$\varphi_m$ on the segment~$[\pi , h_{2k+1})$. 
 If $a$ is large enough, then the function 
 $\varphi$ possesses the same property for a finite 
 sequence~$\pi = h_0'\, < \, h_1' \, < \ldots < h_{2n}'$. 
 Moreover, 
 $n$ can be arbitrarily large, provided~$a$ is big enough. 
Now it remains to set~$\tau_j = h_j'$.

{\hfill $\Box$}
\smallskip 

\bigskip

\begin{center}
\large{\textbf{6. Deciding stability 
and computing the Lyapunov exponent}}
\end{center}
\bigskip

We present an algorithm which for every system 
$S= \{\cA, m, M\}$ with a finite set~$\cA$ decides the stability with some approximation rate. The algorithm decides between two cases: 
either $\sigma(S) > -\delta$ or $\sigma(S) <  \varepsilon$, where 
$\delta > 0$ is a chosen precision and $\varepsilon > 0$ is 
an upper bound that depends on the parameters of the algorithm and that can  also be made 
arbitrarily small. Thus, the algorithm solves the stability problem  
approximately with an arbitrary precision. Then we can compute the Lyapunov exponent 
$\sigma(S)$ by the double division. 

To avoid technical difficulties, we assume that our system is irreducible 
i.e., its matrices do not share a common nontrivial invariant subspace. 
This assumption is made without loss of generality since the stability problem can always 
be reduced to this case. Namely, every reducible system 
is split into several systems of smaller dimension and its Lyapunov exponent 
is equal to the maximal Lyapunov exponents of those systems, see~\cite{B, CGP} for details. 

\medskip

\bigskip

\begin{center}
\textbf{6.1. The main idea and the structure of the algorithm}
\end{center}
\bigskip

We exploit the idea of the invariant polytope algorithm from~\cite{GP13}, 
when a polyhedral Lyapunov norm is constructed  iteratively. 
The theoretical base of the algorithm is provided by three our results: 
\smallskip 

Theorem~\ref{th.20} -- the existence of the Lyapunov norm for stable systems; 

Theorem~\ref{th.40} -- the improved Berger-Wang formula.  

Proposition~\ref{p.80} -- the estimate of the precision of discretization. 
\smallskip 

The algorithm finds the Lyapunov multinorm for the discretized system, 
when the switches on the interval $[m, M]$ are allowed only at the 
points of the uniform partition:~$m + k\tau, k = 0, \ldots , N$, where $\tau = \frac{M-m}{N}$. 
Thus, after the operator~$B_j$ we apply some power of the operator~$e^{\, \tau A_j}$
in the space~$L_j$. Since  the discretized system has a smaller set of trajectories, it follows that its Lyapunov exponent gives a lowed bound for~$\sigma(S)$. On the other hand, 
Proposition~\ref{p.80} proved below in Section~9 gives an inverse estimate: 
the distance between the Lyapunov exponents of these two systems does not exceed 
$\frac{C}{N^2}$, where $C$ is an effective constant. 

The stability of the discretized system is verified by  its Lyapunov multinorm. 
We construct a polyhedral multinorm which is defined in each space $L_j$
as the Minkowski functional of a polytope~$P_j$. The set of polytopes~$\{P_j\}_{j=1}^n$ is  computed iteratively. At the first iteration we take some index~$i$
and a nonzero point~$\bx_0\in L_i$. The further construction is by induction. 
If after the $k$th iteration we have polytopes~$P_j^{(k)} \in L_j, \, j=1, \ldots , n$, 
then for each vertex~$\bv$ of~$P_j^{(k)}$ and for each $s=1, \ldots , n$,
we consider $N$ points in the space~$L_s$: 
$ e^{ r\tau A_s } B_s\, \bv, \, r=0, \ldots , N $. Those points 
split the trajectory of the ODE $\dot \bx(t) \, = \, A_j\bx(t), \, t \in [m, M], \, 
\bx(m)\, = \, B_s\, \bv$ into $N$ arcs. If some of those points are inside the 
polytope~$P_j^{(k)}$, we call them {\em dead} and remove those 
points. All other points are {\em alive} and are added to the set of vertices of~$P_j^{(k)}$. 
Thus, the polytopes~$P_j^{(k)}$ grow with each iteration. Simultaneously 
we compute the spectral radii of the matrix products~$\Pi(\bv)$ corresponding to each alive vertex and get a lower bound for $\sigma(S)$. We consider only allowed products, which 
belong to~$\cE_0$.

The algorithm halts when either the lower bound becomes positive 
(in which case the system is unstable) or no alive vertices appear in a current iteration. 
After  this we get the bounds for $\sigma(S)$ which are provided by  
Theorem~\ref{th.40} (the lower bound) and   
by Proposition~\ref{p.80} (the upper bound). The distance between bounds can be made arbitrarily small. The algorithm always halts, for every 
restricted system (Theorem~\ref{th.50}).

\bigskip

\begin{center}
\textbf{6.2. The  algorithm}
\end{center}
\bigskip

\bigskip 

\noindent \textbf{Preliminary notation and remarks}. 
For a given~$N$, we set 
$$
\tau \ = \ \frac{M-m}{N}, \qquad  \cA_{\tau}\  = \ 
\Bigl\{  e^{s\tau A_j}B_j, \, s = 0, \ldots , N ; \ j = 1, \ldots , n\Bigr\}.
$$ 
If $\cA_{\tau}$ is reducible for all~$N$, then 
the family~$\bigl\{ e^{t A_j}B_j, \ t \in [0, M-m] \bigr\}$ is reducible, 
which is impossible by Lemma~\ref{l.20}. Hence, there exists $N$ such 
$\cA_{\tau}$ is irreducible. In the algorithm we use the notation for the~symmetrized  convex 
  hull~${\rm co}_s\, X \, = \, {\rm co}\, \bigl\{X, -X \bigr\}$;  by~$\Pi(a,b)$ 
  we denote the solution of the matrix equation~$\Pi' = A(t)\Pi, \
  \Pi(a) = I$ and call it {\em product}. Thus, $\Pi(0,T)\, = \, \Pi(T)$. 
\bigskip 

\textbf{Algorithm~1}. 
\medskip 

We have an irreducible  system~$S = \{\cA, m, M\}$ with $\cA = \{A_1, \ldots , A_n\}$. 
For each~$j$, we use notation $L_j$ for the space~$\re^d$ as the 
domain of the operator~$A_j$. 
 Choose the precision~$\delta > 0$.  

\medskip 

 \textbf{Initialization}. 
We fix natural~$N \ge 1$ and denote~$\tau = \frac{M-m}{N}$. 
Verify that $\cA_{\tau}$ is irreducible. 
Choose arbitrary~$i \in \{1, \ldots , n\}$ and a 
starting point~$\bx(0) \in L_i\setminus \{0\}$. 
Define two sets~$\cV_i^{(1)} = \cR_i^{(1)} \, = \, \{\bx(0)\}$; 
for all~$j\ne i$, we set~$\cV_j^{(1)} = \cR_j^{(1)} \, = \, \emptyset$.  

\bigskip

 \textbf{The main loop}. 
$k$th iteration.  
\bigskip 

{\tt In each space~$\, L_j, \ 1, \ldots , n$},  
we have a finite set~$\cV_j^{(k)} \subset L_j$ and its subset~$\cR_j^{(k)}$. 
Define~$\cV_j^{(k+1)}\, = \, \cV_j^{(k)}, \ \cR_j^{(k+1)}\, = \, \emptyset$. 
\medskip 

{\tt For every $j= 1, \ldots , n$}, we take consequently all
$q \in \{1, \ldots, n\}\setminus \{j\}$. 
\medskip 

{\tt For each $q$}, we take consequently all points of the 
set~$\cR_q^{(k)}$. 
\medskip 

{\tt For each point~$\bx \in \cR_q^{(k)}$}, 
we construct $N+1$ points in the space~$L_j$: 
  $\bx_s \, = \, e^{\, s \tau A_j}B_j\bx, \ s = 0, \ldots , N$. 
  For each of them, we check whether $\bx_s$ belongs to the interior of the symmetrized  convex 
  hull~${\rm co}_s\, \cV_j^{(k+1)}$. This is done by the auxiliary Problem 1 (problem~(\ref{eq.LP})). 
  If so, then $\bx_s$ is {\em dead} and  we leave~$\cR_j^{(k+1)}$ and~$\cV_j^{(k+1)}$ as they are. 
  If not, then the point~$\bx_s$ {\em survived} and we add it 
  to both those sets~$\cR_j^{(k+1)}$ and~$\cV_j^{(k+1)}$. 
  Then we consider the trajectory of the system leading from the starting point~$\bx(0)$
  to~$\bx_s$. Denote by~$t$ the time of the point~$\bx_s$, i.e., $\bx(t)\, = \, \bx_s$. 
  Denote by~$t_1, \ldots , t_{\ell}$ the time moments 
  of switches form the space~$L_j$ to other spaces
  (or, the same, from the matrix~$A_j$ to other matrices). Then we set
  $\mu$ to be equal to~$\max\, \bigl\{\mu \, , \, \frac{1}{t-t_i}\sigma\bigl(\Pi(t_i, t)\bigr), \ 
  i = 1, \ldots , \ell \bigr\}$ (if there are no such switches, then we leave~$\mu$ as it is).  
  
  
  \medskip

{\tt If  $s<N$}, then  
  we pass to the next point~$\bx \in \cR_q^{(k)}$. 
 \smallskip  
 
  {\tt If all points of~$\cR_q^{(k)}$
  are exhausted}, then we pass from~$q$ to $q+1$, provided 
  $q+1\ne j$ (otherwise we pass to~$q+2$). 
  \smallskip  
  
  {\tt If all~$q \in \{1, \ldots , n\}\setminus \{j\}$
  are exhausted}, then we pass from~$j$ to $j+1$. 
  \smallskip  
  
  {\tt If $j=n$}, then the $k$th iteration is completed. 
  
\bigskip 
  
\noindent  \textbf{Termination}. There are two possible cases: 
 \medskip 
 
\noindent    \textbf{Case 1}. If $\mu > -\delta$ at some moment of the algorithm, then the algorithm halts. We have~$\sigma(S)\ge -\delta$. 
  \smallskip 
  
\noindent    \textbf{Case 2}. 
  If~$\cR_j^{(k+1)} \, = \, \emptyset$ for all~$j$, i.e., the $k$th iteration does not produce any new point, then the algorithm halts. Denote~$P_j \, = \, {\rm co}_s\, \cV_j^{(k)}$
  and find the value
  $$
  \bigl\|A_j^2\bigr\|_{P_j}\quad = \quad \min \ \Bigl\{\, \lambda > 0\quad  \Bigl| \quad  
  A_j^2 P_j \subset  \lambda P_j \ \Bigr\}\, , 
  $$ 
  which is the operator norm of~$A_j^2$ induced by the Minkowski norm of the polytope~$P_j$. 
 Let~$\|\cA^2\|_P \, = \, \max\limits_{j=1, \ldots , n} 
 \bigl\|A_j^2\bigr\|_{P_j}$.  This value is computed by the auxiliary Problem 3 
 given below. 
 Then~$\sigma(S) < \nu$, where  
 \begin{equation}\label{eq.est1}
\nu  \ = \ -\, \frac1m\, \ln\, \left(1\, - \, 
\frac{(M-m)^2\|\cA^2\|_P}{8N^2} \right)\, . 
 \end{equation}
   \medskip

 \textbf{End of the Algorithm~1}
  \bigskip

\noindent \textbf{The auxiliary problems}. Given a finite set~$\cV = \{\bx_i\}_{i=1}^{\ell}$
and a point~$\bx_0$. Consider the polytope~$P = {\rm co}_s\, \cV$. 
\medskip 

\noindent  \textbf{1)} Decide whether~$\bx_0 \in {\rm int}\, P$; 
\smallskip 

\textbf{Solution.} Solve the linear programming (LP) problem: 
\begin{equation}\label{eq.LP}
\left\{
\begin{array}{l}
t_0 \  \to \  \max\\
t_0 \bx_0 \ = \ \sum_{j=1}^{\ell} (t_i - s_i)\bx_i\\
\sum_{j=1}^{\ell} t_i \ = \ \sum_{j=1}^{\ell}  s_i \ = \ 1\\
t_i \ge 0, \ s_i \ge 0, \qquad i=1, \ldots , \ell\, . 
\end{array}
\right. 
\end{equation}
Then $\bx_0 \in {\rm int}\, P\ $ iff $\ t_0 > 1$. 
\medskip
 
\noindent  \textbf{2)}.  
Find $\|\bx_{0}\|_P$.

\smallskip  \textbf{Solution.} ${\|\bx_{0}\|_P \, = \, 1/t_0}$, where $t_0$
is the solution of the LP problem~(\ref{eq.LP}). 

\smallskip
 
\noindent  \textbf{3)}.  
For a given~$d\times d$ matrix~$A$, find $\|A\|_P$.
\smallskip 

 \textbf{Solution.} 
For every~$j= 1, \ldots , \ell$, set~$\bx_0 = A\bx_j$, solve the LP problem~(\ref{eq.LP})
and denote~${t_0 = t_0^{(j)}}$. Then~$\|A\|_P \, = \, \max_{j=1}^{\ell} \|A\bx_j\|_P\, = \, 
\max_{j=1}^{\ell} 1/t_0^{(j)}$. 

\bigskip

The main result of this section states that Algorithm~1 always terminates within finite time. 

\begin{theorem}\label{th.50} 
For every initial data, Algorithm~1 halts within finite time. 
If it halts in Case~1, then~$\sigma(S) > -\delta$.
Otherwise, it halts in Case~2 and $\sigma(S)< \nu$, where $\nu$ is defined in~(\ref{eq.est1}), and the multinorm~$\|\cdot \|_{P}$ is Lyapunov for the system~$\cA - \nu I$.  
\end{theorem}
The proof is in Section~9.
\begin{remark}\label{r.50}
{\em The parameter~$\delta$ can be chosen arbitrarily small, which could 
allow Algorithm~1 to compute the Lyapunov exponent with an arbitrary precision. 
However, there is also parameter~$\nu$ which increases the distance 
between upper and lower bounds. This is unavoidable in the framework of Algorithm~1 
since it constructs the polyhedral Lyapunov norm for the discretized system.  
The parameter~$\nu$ is the precision of approximation by the discretization. 
It is asymptotically equivalent to $ 
\frac{1}{8N^2}\, \frac{(M-m)^2\|\cA\|_P}{m}$ as $N \to \infty$.
Increasing~$N$ we improve the precision, but increase the number of vertices, 
which slows down the algorithm.

}
\end{remark}

\noindent \textbf{Computing the Lyapunov exponent~$\sigma(S)$}. 
By Theorem~\ref{th.50},  Algorithm~1 decides within finite time 
between two cases: either~$\sigma(S) > -\delta$ or $\sigma(S)< \nu$.
Therefore, introducing a parameter~$\alpha$ and applying Algorithm~1
to the system~$\cA - \alpha I$, one can realize the double 
division in~$\alpha$. The initial interval for $\alpha$
can be chosen, for example, as $\, \bigl[ \, \max_{j=1}^n \sigma(A_j)\ , \  
\hat \sigma(S) \bigr]$, where $\hat \sigma (\cA)$ is an arbitrary upper bound for the unrestricted system with  matrices~$\cA$, for example, it can be the CQLF estimate~\cite{LHM}
or polyhedral estimates~\cite{BM, GLP17}. 
Let $L$ be the length of the initial interval. Then after $k\ge 1$ iterations  Algorithm~1 localizes the value of $\sigma (\cA)$ on an interval of 
length~$2^{-k}L + \delta + \nu$. The parameter $\delta$ is under our control
and can be chosen arbitrarily small. For the value 
$\nu$ defined in~(\ref{eq.est1}), we have~$\nu = O\bigl(\frac{1}{N^2} \bigr)$, so it can 
be made small by increasing~$N$. Thus, Algorithm~1 along with the double division and 
with the proper choice of the parameters~$\delta$ and $N$ computes 
the Lyapunov exponent with a given precision. 
\medskip 

\noindent \textbf{Algorithm~1 for positive systems}.  
If the restricted system~$S = \{\cA, m, M\}$ is positive, then Algorithm~1 can be modified, 
which in practice gives much more efficiency. Let us recall that the system is called 
positive if all the matrices from~$\cA$ are Metzler, i.e., all their off-diagonal elements
are nonnegative. In this case all trajectories of the system starting on the 
positive orthant~$\re^d_+$ stays inside it for all~$t \in [0, +\infty)$. 
See, \cite{FM, GLP17} for more on positive systems.

The modification is the following.  For an arbitrary subset
$K$ of the positive orthant~$\re^d_+$, we consider its 
{\em positive convex hull} ${\rm co}_{p}\, K\, = \, \bigl\{\bx \ge 0,\, 
\exists \, \by \in {\rm co}\, K  \ \mbox{s.t.} \ \by \ge \bx\bigr\}$, 
where the vector inequality $\by \ge \bx$ is coordinatewise. 
If $K$ is a finite set, then ${\rm co}_{p}\, K\, $ is a {\em positive polytope}. 
 In Algorithm~1 we replace the polytopes~$P_j = {\rm co}_{s}\, \cV$ by the 
positive polytopes~$P_j = {\rm co}_{p}\, \cV$. In  Algorithm~1 
we everywhere replace the symmetrized convex hill by the positive 
convex hull, respectively, the LP problem~\ref{eq.LP} is replaced by  
\begin{equation}\label{eq.LP+}
\left\{
\begin{array}{l}
t_0 \  \to \  \max\\
t_0 \bx_0 \ = \ \, \by\, + \, \sum_{j=1}^{\ell} t_i\bx_i,\\
\sum_{j=1}^{\ell} t_i \ = \ 1,\\
\by \ \le \ 0,\\
t_i \ge 0, \ \qquad i=1, \ldots , \ell\, . 
\end{array}
\right. 
\end{equation}
As a rule, this version of Algorithm~1 leaves much less alive vertices of polytopes~$P_j^{(k)}$ (in this case, positive polytopes) and converges much faster, even in high dimensions. 

\medskip

  If Algorithm~1 does not terminate within reasonable time, it can be interrupted, 
  in which case we have the following output: 
 \smallskip 
  
  \textbf{Complement to Algorithm~1.}
  \medskip
Choose~$K$ as the maximal number of iterations.  
   If Algorithm~1 does not halt by the~$K$th iteration, then we interrupt it and make one more iteration.  For every survived point~$\bx_s$ generated in this iteration 
   (i.e., $\bx_s \in L_i, \, \bx_s  \notin {\rm int} \, P^{(K)}_i$).
    Then we denote by~$\gamma$
   the biggest one of those numbers and obtain
 \begin{equation}\label{eq.est2}
\mu \ \le \  \sigma(S) \ \le \ -\, \frac1m\, \ln\, \left(1 \, - \, 
\frac{(M-m)^2\gamma\|\cA^2\|_P}{8N^2} \right)\, . 
 \end{equation}  
 \smallskip   
  
  \bigskip 

\newpage 

\begin{center}
\large{\textbf{7. Numerical results}}
\end{center}
\bigskip

We consider two low-dimensional examples of computing the Lyapunov exponent 
by Algorithm~1. Then we present statistics of numerical experiments 
for higher dimensions for general restricted systems and for positive systems. 
\bigskip 

\begin{center}
\textbf{7.1. Examples}
\end{center}
\bigskip

\begin{ex}\label{ex.10num}
{\em 
	For the system with $\cA = \{A_1 , A_2\}, \ m=1, M =2 $, where 
	\[ 
	A_1 = \left(
	\begin{array}{cc}
		-0.3 & 0.5 \\
		0.2 & -0.4
	\end{array}
	\right),\,
	A_2 = \left(
	\begin{array}{cc}
		-0.6 & 0 \\
		0 & 1
	\end{array}
	\right),  \qquad m = 1, M = 2.
	\]
	 Algorithm~1 
	gives the following results:  
\begin{table}[h!]
	\caption{$ 2 \times 2$ matrices , $ m = 1, M = 2 $.}
	\label{tabular1:timesandtenses}
	\begin{tabular}{|c|c|c|c|c|c|}
		\hline
		$ \text{len} P_1 $ & $ \text{len} P_2 $ & $ N $ & $ \text{\# iterations} $ & $ \sigma $ & $ \text{time} $ 
		\\
		\hline
		$ 28 $ & $ 34  $ & $ 10 $  & $ 16 $ & $ (0.15560 , 0.15685) $ & $ 11.9$
		\\ 
		\hline
	\end{tabular}
\end{table}	
Here~$\text{len} P_i$ denotes the half of the number of verices 
of the polytope~$P_i$ that generates the Lyapunov norm~$f_i, \, i = 1, 2$;
$N$ is the number of points of discretization (so, the discretization step is 
$\tau = \frac{M-m}{N} = 0.1$); the column ``$\sigma$'' shows the interval 
for $\sigma(S)$; ``time'' shows the time (sec) of running Algorithm~1 in a standard laptop. 

We see that in this example Algorithm~1 gives the estimates 
$0.15560 < \sigma(S) < 0.15685$, the relative error is less than $1.3\cdot 10^{-3}$.  
}
\end{ex}
\begin{ex}\label{ex.20num}
{\em 
For the system with $\cA = \{A_1 , A_2\}, \ m=1, M =2.5 $, where 
\[
A_1 = \left(
\begin{array}{cc}
	-0.3 & 0.5 \\
	0.2 & -0.4
\end{array}
\right),\,
A_2 = \left(
\begin{array}{cc}
	-0.6 & 0 \\
	0 & 1
\end{array}
\right),  \qquad m = 1, M = 2.5.
\]
	 Algorithm~1 
	gives the following results:  
\begin{table}[h!]
	\caption{matrices $ 2 \times 2$, $  m = 1, M = 2.5 $.}
	\label{tabular2:timesandtenses}
	\begin{tabular}{|c|c|c|c|c|c|}
		\hline
		$ \text{len} P_1 $ & $ \text{len} P_2 $ & $ N $ & $ \text{\# iterations} $ & $ \sigma $ & $ \text{time} $ 
		\\
		\hline
		$ 10 $ & $ 8 $ & $ 10 $  & $ 10 $  & $ ( 0.61507,  0.62433) $ & $  5.4 $
		\\ 
		\hline
	\end{tabular}
\end{table}	
The interval for~$\sigma$ is  $( 0.61507,  0.62433)$, the relative error is 
less than $0.01$. 

}
\end{ex}

\bigskip

\newpage 

\begin{center}
\textbf{7.2. Results of numerical experiments}
\end{center}
\bigskip

Table~\ref{tabular4:timesandtenses} shows the results of Algorithm~1 
applied to systems with two $d\times d$ matrices with~$[m,M] = [1, 2]$.  
The notation are the same as in Example~\ref{ex.10num}, $d$ is the dimension.

\begin{table}[h!]
	\caption{$d\times d$ matrices  $  m = 1, M = 2 $.}
	\label{tabular4:timesandtenses}
	\begin{tabular}{|c|c|c|c|c|c|c|c|}
		\hline
		$ \text{dim} $ & $ \text{len} P_1 $ & $ \text{len} P_2 $ &  $ N $ & $ \text{\# iterations} $ & $ \sigma $ & $ \text{time} $ 
		\\
		\hline
		$ 3 $ & $ 16 $ & $ 10 $ & $ 10 $  & $ 18 $  & $ (1.15482, 1.16474)  $ & $  26.5 $
		\\ 
		\hline
		$ 4 $ & $ 44 $ & $ 46  $ & $ 10 $  & $ 22 $  & $ (2.52396, 2.98327) $ & $ 128.1  $
		\\ 
		\hline
			$ 5 $ & $ 138 $ & $ 106  $ & $ 10 $  & $ 12 $  & $ (3.00530, 4.03136) $ & $  96.9  $
		\\ 
	    \hline
	
	\end{tabular}
\end{table}	

Table~\ref{tabular5:timesandtenses} shows the results of Algorithm~1 
applied to systems with two $d\times d$ matrices with~$[m,M] = [1, 2.5]$.  

\begin{table}[h!]
	\caption{$d\times d$ matrices, $  m = 1, M = 2.5 $.}
	\label{tabular5:timesandtenses}
	\begin{tabular}{|c|c|c|c|c|c|c|}
		\hline
		$ \text{dim} $ & $ \text{len} P_1 $ & $ \text{len} X_P $ & $ N $ & $ \text{\# iterations} $ & $ \sigma $ & $ \text{time} $ 
		\\
		\hline
		$ 3 $ & $ 10 $ & $ 12 $ & $ 10 $ & $ 12 $  & $ (1.23099,  1.25315)  $ & $  11.2 $
		\\ 
		\hline
		$ 4 $ & $ 60 $ & $ 62  $ &  $ 15 $ &  $ 16 $  & $ ( 2.57561, 3.06883) $ & $ 126.1 $
		\\ 
		\hline
		$ 5 $ & $ 144 $ & $ 140  $ &  $ 15 $ &  $ 8 $  & $ (3.15519, 4.10713) $ & $ 188.9 $
		\\ 
		\hline
		
	\end{tabular}
\end{table}

Tables~\ref{tabular6:timesandtenses} and~\ref{tabular7:timesandtenses}
  show the results for positive systems. We see that in this case we have 
 quite sharp computations even for relatively hight dimensions, up to~$35$. 

 \begin{table}[h!]
	\caption{nonnegative matrices,  $  m = 1, M = 2$.}
	\label{tabular6:timesandtenses}
	\begin{tabular}{|c|c|c|c|c|c|c|}
		\hline
		$ \text{dim} $ & $ \text{len} P_1 $ & $ \text{len} P_2 $ &  $ N $ &$ \text{\# iterations} $ & $ \sigma $ & $ \text{time} $ 
		\\
		\hline
		$ 3 $ & $ 1 $ & $ 2 $ & $ 10 $ & $ 6 $  & $ 9.02606  $ & $  46.5 $
		\\ 
		\hline
		$ 5 $ & $ 1 $ & $ 2 $ & $ 10 $ &  $ 6 $  & $  13.62821 $ & $ 35.3  $
		\\ 
		\hline
		$ 10 $ & $ 1 $ & $ 1  $ & $ 10 $ &  $ 4 $  & $ (24.70551, 24.70731) $ & $  5.3$
		\\ 
		\hline
		$ 20 $ & $ 1$ & $ 1  $ & $ 10 $ &  $ 4 $  & $ (50.07293, 50.07355) $ & $  4.5 $
		\\ 
		\hline
		$ 35 $ & $ 1 $ & $ 1  $ & $ 10 $ &  $ 4 $  & $ (87.07711, 87.0775) $ & $  7.8$
		\\ 
		\hline

	\end{tabular}
\end{table}	

 \begin{table}[h!]
	\caption{nonnegative matrices,  $  m = 1, M = 2.5 $.}
	\label{tabular7:timesandtenses}
	\begin{tabular}{|c|c|c|c|c|c|c|}
		\hline
		$ \text{dim} $ & $ \text{len} P_1 $ & $ \text{len} P_2 $ &  $ N $ &$ \text{\# iterations} $ & $ \sigma $ & $ \text{time} $ 
		\\
		\hline
		$ 3 $ & $ 1 $ & $ 2 $ & $ 10 $ & $ 6 $  & $  9.17010  $ & $ 29.1 $
		\\ 
		\hline
		$ 5 $ & $ 1 $ & $ 2 $ & $ 10 $ &  $ 6 $  & $ 13.66565 $ & $ 8.8  $
		\\ 
		\hline
		$ 10 $ & $ 1 $ & $ 1  $ & $ 10 $ &  $ 4 $  & $ (24.87671, 24.87968) $ & $  5.3 $
		\\ 
		\hline
		$ 20 $ & $ 1$ & $ 1  $ & $ 10 $ &  $ 4 $  & $ (50.12228, 50.12387 $ & $  4.0 $
		\\ 
		\hline
		$ 35 $ & $ 1 $ & $ 1  $ & $ 10 $ &  $ 4 $  & $ (87.11728, 87.11831) $ & $  5.3 $
		\\ 
		\hline

	\end{tabular}
\end{table}

\bigskip 

\newpage

\begin{center}
\large{\textbf{8. When can the upper bound be cancelled?}}
\end{center}
\bigskip

In this section we deal with the general case of mode-dependent 
segments~$[m_j, M_j]$. 
Algorithm~1 may suffer when the maximal switching time~$M_j$ is big. 
Indeed, by Theorem~\ref{th.50}, the upper bound for the Lyapunov exponent given by~(\ref{eq.est1}) requires a large number~$N$ of discretization segments 
 to get a satisfactory precision. 
This significantly increases the number of vertices of the polytopes and, therefore,  the running time.
The question arises whether the upper bound~$M_j$ can be reduced without 
changing the stability/instability of the system?  Reducing  some of those bounds simplifies the implementation of Algorithm~1. 

 On the other hand, in case~$M_j = +\infty$, when there is no 
upper bound for the switching time of~$A_j$, the stability can be decided by many methods 
 known from the literature~\cite{CGPS21}. So, the opposite problem 
is when the upper bound can be cancelled without influencing the stability?

 \smallskip 
 
  \noindent \textbf{Problem 1.} {\em For a given $j$, under what  conditions
  the upper time limit  $M_j$ can be reduced (and to which point) without changing the stability/instability of the system? }
 
  \medskip 
 
 \noindent \textbf{Problem 2.} {\em For which $M_j$ the stability of the system 
 does not change 
 after replacing~$M_j\, $ by $\, +\infty$ ?  }
 
  \smallskip

 We attack these problem by considering first a system~$\dot \bx\, = \, A\bx, \bx(0) = \bx_0$
with a constant matrix~$A$,  without switches. This system is stable when~$A$ is stable, or a Hurwitz matrix, i.e., when $\sigma(A) < 0$.  

In what follows we assume that 
 $\bx_0$ is a {\em generic point}, i.e., it does not belong to any proper invariant subspace of~$A$. For an arbitrary segment~$[t_1, t_2] \subset \re_+$, we   
 denote~$\Gamma(t_1, t_2)\, = \,\bigl\{\bx(t) \ | \ 
t \in [t_1, t_2]\bigr\}$ and $G(t_1, t_2)\, = \, {\rm co}_s \, \Gamma(t_1, t_2)$.  
We also denote by  $\Gamma = \Gamma(0, +\infty)$ the entire trajectory and, respectively, 
$G = G(0, +\infty)$.  
 It is seen easily that $G(t_1+h, t_2+h)\, = \, 
e^{hA}\, G(t_1, t_2)$. 
 \begin{defi}\label{d.80}
 A point~$T>0$ is called a {\em cut tail point} for 
 the stable system~$\dot \bx\, = \, A\bx, \ \bx(0) = \bx_0$, if  
 for every $t> T$, the point $\bx(t)$  belongs to the interior of the set~$G(0, T)$. 
 \end{defi}   

 \begin{prop}\label{p.100}
For every stable system~$\dot \bx\, = \, A\bx, \ \bx(0) = \bx_0$, the set of cut tail points is nonempty. 
 \end{prop}  
 {\tt Proof.}  Since the point $\bx_0$ is generic, it follows that every arc of the trajectory is not contained in a hyperplane, therefore,  for each segment~$[t_1, t_2]$, 
the set $G(t_1, t_2)$
is full-dimensional and hence,  
it contains a ball centered at the origin. 
On the other hand, the trajectory converges to zero as~$t\to \infty$, hence, 
say, $G(0,1)$ contains $\bx(t)$ for all~$t$ bigger than some $T>1$. 
Since $G(0,1) \subset G(0,T)$, we see that~$T$ is a cut tail point.  
{\hfill $\Box$}
\medskip

The set  of cut tail points is obviously closed, hence, there exists
 the minimal point, which will be denoted by~$T_{cut}$. Clearly, 
$\bx(T_{cut})$ is on the boundary of~$G$. 
If two trajectories 
start at generic  points~$\bx_0, \bx_0'$, then 
they are affinely similar. 
Hence {\em cut tail points do not depend on~$\bx_0$} 
and are functions of the matrix~$A$. Thus, speaking about 
the cut tail points we may not specify~$\bx_0$.

\begin{prop}\label{p.90}
Suppose~$A$ is a stable matrix; then the arc~$\Gamma (0, T_{cut})$
lies on the boundary of~$G$.  
\end{prop}
{\tt Proof.} If, on the contrary, there exists~$t< T_{cut}$ such 
that~$\bx(t) \in {\rm int} \, G$, then, since the ``tail'' $\{\bx(s)\ | \ s> T_{cut}\}$
of the trajectory~$\Gamma$
lies in the interior of~$G$, 
it follows that~$\bx(t) \in {\rm int} \, G(0, T_{cut})$. 
Hence,
$$
\bx(T_{cut}) \ = \,  e^{(T_{cut} - t)A}\, \bx(t)\ 
\in \ {\rm int} \, e^{(T_{cut} - t)A}\, G(0, T_{cut})\ = \ {\rm int} \, 
G\, (T_{cut} - t\, , \, 2T_{cut}-t)\ \subset \ {\rm int}\, G.
$$
Thus, $\bx(T_{cut}) \in {\rm int}\, G$, which is a contradiction. 

{\hfill $\Box$}
\medskip 

We see that the point  $\bx(T_{cut})$ separates the part of trajectory that lies on the 
boundary form that in the interior: $\bx(t) \in \partial G$, for all~$t\le T_{cut}$, and 
$\bx(t) \in {\rm int}\, G$, for all~$t > T_{cut}$. 
\begin{cor}\label{c.34}
The set of cut tail points is precisely the half-line
$[T_{cut}, +\infty)$. 
\end{cor}

Now turn back to the switching systems. Consider  a stable system  $S = \{\cA, \bm, \cM\}$. 
Let us omit the upper limit~$M_j$ for some of the matrices~$A_j$, i.e.,
we allow the switching time for those matrices to run over~$[m_j, +\infty)$. 
The new system  has a bigger set of of trajectories, so 
its Lypunov exponent is at least not smaller. In particular, 
 the Lypunov exponent can become positive, in which case the system loses stability.
 We are going to show that if~$T = M_j-m_j$ is a cut tail point for the matrix~$A_j$, 
 then the system stays stable and, moreover, have the same Lyapunov multifunction.   
 
\begin{theorem}\label{th.60}
Let $S=\{\cA, \bm, \cM\}$ be a restricted system with all the matrices~$A_j$
being stable. Let~$J$ be an arbitrary subset of the index set of~$\cA$. 
 If for each $j\in J$, we have 
$T_{cut}(A_j) \le M_j-m_j$, then 
replacing each~$M_j$ either by $+\infty$ or by~$m_j+T_{cut}(A_j)$
changes neither the stability of~$S$ nor  
the  Lyapunov multinorms. 
\end{theorem}
Thus, if the number~$M_j-m_j$ is a cut tail point for some of the matrices~$A_j$, then 
one either omit this upper bound or reduce it to~$m_j + T_{cut}(A_j)$
without changing the stability/instability. 
 \begin{cor}\label{c.50}
If for a stable system~$S=\{\cA, m, M\}$ with stable matrices, we have~$M - m\, \ge \, \max_{j=1, \ldots , n}
T_{cut}(A_j)$, then the system $\{\cA, m, +\infty\}$ is also stable and have the same 
Lyapunov multinorms. 
\end{cor}
{\tt Proof of Theorem~\ref{th.60}}. If~$S$ is stable, then it possesses a 
Lyapunov multinorm~$f$.  
Then for each~$i \in J$ and for every~$k \ne j$ and for every~$\bx \in \cS_k$, 
 the curve 
$\gamma \, = \, \bigl\{ e^{t A_{i}}\, B_{i}\, \bx \, , \ t \in [0, M_i-m_i]\, \bigr\}$,  
is strictly inside the ball~$\cB_{i}$
(Proposition~\ref{p.20}). If $M_j - m_j \, \ge \, T_{cut}(A_{i})$, then~$M_i-m_i$  is a cut tail point for~$A_i$. 
 Therefore, the entire curve~$e^{t A_{i}}\, B_{i}\, \bx \, , \ t \in [0, +\infty)$
is contained in ${\rm co}_s \, \Gamma$ and  hence, it lies inside~$\cB_i$.   
Thus, if $i \in J$, then $f$  is a Lyapunov multinorm for the system 
with the omitted upper restriction for~$A_{i}$. Similarly we show that 
if after the reduction of~$M_i$ to~$T_{cut}(A_{i})$, the system 
is stable, then the original system is stable as well.

{\hfill $\Box$}
\medskip 

Now the problem becomes how to decide if a given~$T$ is a cut tail point for 
a given stable matrix~$A$ and how to find~$T_{cut}(A)$?  Let us denote by 
$\cP_{A}$ the linear span of the functions $f_{\bx}(t) = e^{tA}\bx, \  \bx \in \re^d$. 
This is the  space of quasipolynomials 
which are linear combinations of functions~$t^ke^{\, \alpha t}
\sin \beta t$ and $t^ke^{\, \alpha t}\, \cos \beta t\, $, 
where $\alpha + i \beta$ is an eigenvalue of~$A$ and $k = 1, \ldots , r-1$, where $r$ is the size of the 
largest Jordan bock of~$A$ corresponding to that eigenvalue.  The dimension of~$\cP_A$
is equal to the degree of the minimal annihilating polynomial of the matrix~$A$. 
\begin{theorem}\label{th.70}
Let~$A$ be a stable matrix; then a number~$T>0$ is a cut tail point if and only if 
the value of the 
convex extremal problem 
\begin{equation}\label{eq.extr-cut}
\left\{
\begin{array}{l}
p(T) \to \max \\
\|\bp\|_{C(\re_+)} \ \le \ 1\\
p \in \cP_A
\end{array}
\right. 
\end{equation}
is smaller than~$1$.  
\end{theorem}
{\tt Proof}. By Proposition~\ref{p.90}, the inequality~$T> T_{cut}$ is equivalent to that 
$\bx(T) \in {\rm int}\, G$. This means that the point $\bx(T)$
cannot be separated from~$G$ by a linear functional, i.e.,  
for every nonzero vector~$\bp\in \re^d$, we have 
$\bigl(\bp , \bx(T)\bigr) \, <   \, \sup\limits_{\by  \in G} \, (\bp , \by)$. 
Since~$G$ is a convex hull of points $\, \pm \, \bx(t)\, , \, t \ge 0$, 
it follows that $\sup\limits_{\by  \in G} \, (\bp , \by)\, = \, 
\sup\limits_{t \ge 0} \, \bigl|(\bp , \bx(t))\bigr|$. Thus,  
$$
\bigl(\bp , \bx(T)\bigr) \ <   \ \sup\limits_{t \in \re_+}\bigl| (\bp , \bx(t)) \bigr|.
$$ 
On the other hand, $(\bp , \bx(t)) = p(t)$, where $p \in \cP_{A}$
is the quasipolynomial with the vector of coefficients~$\bp$.  
We obtain $p(T) < \|p\|_{C(\re_+)}$. 
 Consequently, $p(T)< 1$
for every quasipolynomial from the unit ball~$\|p\|_{C(\re^+)} \le 1$.
Now by the compactness argument we conclude that the value of the 
problem~(\ref{eq.extr-cut}) is smaller than one.

{\hfill $\Box$}
\medskip


 The problem~(\ref{eq.extr-cut}) is a convex problem on~$\re^d$ and can be solved by the 
 convex optimisation  methods.  For $2\times 2$ matrices, the parameter $T_{cut}$ can be evaluated in the explicit form. 
 \medskip

 \noindent \textbf{The two-dimensional systems}. 

\smallskip 

\noindent {\tt Case 1}. The eigenvalues of $A$ are real. 
For the sake of simplicity we assume they are different, denote them by $\alpha_1, \alpha_2$. 
 The trajectory~$\Gamma =  \{\bx(t)\, | \, t \ge 0\}$ in the basis of eigenvectors of~$A$ 
has the equation~$(x_1(t),  x_2(t)) \, = \, 
\bigl( e^{\alpha_1 t} \, , \,   e^{\alpha_2 t}\bigr)\, , \ 
t \in \re_+$. 
Since $A$ is stable, both~$\alpha_1, \alpha_2$ are negative 
and~$\bx(0) = (1, 1), \, \bx(+\infty) = (0, 0)$. 
Let~$\ba \in \Gamma$ be the point where the tangent line 
drawn from the point~$(-1, -1)$ touches~$\Gamma$. 
Then $G \, = \, {\rm co}_s \, \Gamma$ is bounded by the line segment 
from $(-1, -1)$ to $\ba$ and by the arc of~$\Gamma$ from $\ba$ to 
$(1, 1)$, then reflected about the origin. 
Hence, if $\ba = \bx(T)$, then $\bx(t) \in {\rm int}\, G$ for all~$t > T$, and 
so $T = T_{cut}$. We have~$\, \ba \, + \, s \, \dot \bx(T_{cut}) \, = 
\, (-1,-1)$, where $s> 0$ is some number. 
Writing this equality coordinatewise, we obtain the equation for~$t = T_{cut}$: 
$$
\left\{
\begin{array}{lcl}
e^{\alpha_1 t}\ + \ s\, \alpha_1\, e^{\alpha_1 t}& = & - 1\\
e^{\alpha_2 t}\ + \ s\, \alpha_2\, e^{\alpha_2 t}& = & - 1
\end{array}
\right. 
$$
which becomes after simplification 
$\frac{1 + e^{- \alpha_1 t}}{\alpha_1}\, = \, \frac{1 + e^{- \alpha_2 t}}{\alpha_2}$. 
The unique solution is~$T_{cut}$. 
\medskip

 If the eigenvalues of $A$ are non-real, they are 
 complex conjugate~$\alpha \pm i \beta$ with $\alpha < 0$. 
 The trajectory~$\Gamma$ in a suitable basis 
 has the equation~$(x_1(t), x_2(t)) \,  = \, 
 e^{\alpha t}\, (\cos \beta t \, , \, \sin \beta t)
 \, , \ 
t \in \re_+$. The trajectory $\Gamma$
 goes from the point~$\bx(0) = (1, 0)$ to zero making infinitely  many rotations. 
Taking the point of tangency~$\ba$ of $\Gamma$ with the line 
going from the point~$(-1, 0)$ and arguing as above we obtain 
$$
\left\{
\begin{array}{lcr}
e^{\alpha t} \cos \beta t \ + \ s\, e^{\alpha t}\bigl(\alpha \cos \beta t \, - \, 
\beta \sin \beta t \bigr)& = & - 1\\
e^{\alpha t} \sin  \beta t \ + \ s\, e^{\alpha t}\bigl(\alpha \sin \beta t \, + \, 
\beta \cos \beta t \bigr)& = & 0
\end{array}
\right. 
$$
from which it follows  $\alpha \sin \beta t \, + \, \beta \cos \beta t\, + 
\, \beta \, e^{\alpha t} \, = \, 0$. The unique solution of this equation is $t = T_{cut}$.

 \begin{center}
\large{\textbf{9. Proofs of the fundamental theorems}}
\end{center}
\bigskip 

 In this section we give proofs of Theorems~\ref{th.20},~\ref{th.30},~\ref{th.40}, and \ref{th.50}. 
 
 \bigskip

\begin{center}
\textbf{9.1 Preliminary results. Discrete systems on graphs}
\end{center}
\bigskip

We begin with proof of Theorem~\ref{th.30} on the existence of invariant multinorm 
for irreducible systems. Then we prove Theorems~\ref{th.20} and~\ref{th.20}. 
The proof of Theorem~\ref{th.50} is the most difficult, we give it is a separate subsection 
9.3. 

The proof of Theorem~\ref{th.30} uses some results from 
the discrete dynamical systems on graphs. The theory of such 
systems was developed recently in~\cite{CGP,K,PhJ1}. We begin with basic definitions.

 We have a directed multigraph $G$  
with $m$ vertices $g_1, \ldots , g_n$. Sometimes, the vertices will be denoted by their numbers.
To each  vertex~$i$ we associate a linear space $L_i$ of dimension $d_i \ge 1$.
If the converse is not stated, we assume $d_i \ge 1$.  The set of spaces
$L_1, \ldots , L_n$ is denoted by $\cL$. For each vertices $i, j \in G$
(possibly coinciding),
there is a set $\ell_{ji}$ of edges
from $i$ to $j$. Each edge from $\ell_{ji}$ is identified with a linear operator $A_{ji}: \, L_i \to L_j$.
The family of those operators (or edges) is denoted by $\cA_{ji}$. If  $\ell_{ji} = \emptyset$, then $\cA_{ji} = \emptyset$. Thus, we have a family of spaces~$\cL$
and a family of operators/edges  $\cA = \cup_{i, j}\cA_{ji}$ that act between these spaces according to the
multigraph $G$. This {\em triplet $\, \xi = (G, \cL, \cA)$}  of the multigraph, spaces, and operators  will be referred to as {\em a system on graph}.
A path $\alpha$ on the multigraph $G$ is a sequence of connected subsequent edges, its length (number of edges)
is denoted by $|\alpha|$. The length of the empty path is zero. To every path $\alpha$ along
vertices $i_1\rightarrow i_2 \rightarrow \cdots \rightarrow i_{k+1}$ that consists of edges (operators)
$A_{i_{s+1}i_s} \in \cA_{i_{s+1}i_s}, \, s = 1, \ldots , k$, we associate the corresponding product (composition) of operators $\Pi_{\alpha} = A_{i_{k+1}i_k}\cdots A_{i_2i_1}$.
Note that $|\alpha| = k$. Let us emphasize that a path is not a sequence of vertices but edges. If $G$ is a graph, then any path is uniquely defined by the sequence of its vertices, if
$G$ is a multigraph, then there may be many paths corresponding to the same sequence of vertices.
If the path is closed ($i_1 = i_{k+1}$), then $\Pi_{\alpha}$ maps the space
$L_{i_1}$ to itself. In this case $\Pi_{\alpha}$ is given by a square matrix, and possess eigenvalues, eigenvectors and the spectral radius $\rho(\Pi_{\alpha})$, which is the maximal modulus of its eigenvalues. The set of all
closed paths will be denoted by $\cC (G)$. For an arbitrary $\alpha \in \cC(G)$ we denote by
$\alpha^k = \alpha \ldots \alpha$ the $k$th power of $\alpha$. A closed path is called {\em simple} if it is not a power
of a shorter path.

A multinorm on a  multigraph is introduced in the same way as 
 in Definition~\ref{d.10}.

For a given $\bx_0 \in L_i$ and for an infinite path $\alpha$ starting at the vertex $i$, we
 consider the {\em trajectory} $\{\bx_k\}_{k \ge 0}$ of the system along this path. Here
 $x_{k} = \Pi_{\, \alpha_{k}}\, \bx_0$, where $\alpha_{k}$ is a prefix of~$\alpha$ of length~$k$.
The concepts of stability and of Lyapunov exponent  are defined in a standard way.
The invariant norm is also defined similarly. For the sake of simplicity we will 
reduce the definition to the case~$\sigma(\xi) = 0$. 
\begin{defi}\label{d.60}
Let us have a discrete system on graph~$\xi$ with $\sigma = 0$. 
A multinorm $\|\cdots \|$ is invariant for~$\xi$
 if for every $i$ and $\bx \in L_i$, we have
\begin{equation}\label{eq.invar}
\max_{A_{ji} \in \cA_{ji}, \, j=1, \ldots , n} \, \|A_{ji}\bx\|_j  \ = \, \|x\|_i\, .
\end{equation}
\end{defi}
Similarly to Theorem~\ref{th.30}, the existence of the invariant norm 
holds under the irreducibility assumption. However, the definition 
of irreducible systems on graphs is more complicated, it requires 
not one invariant subspace but a collection of subspaces.  
 A triplet 
$\xi' = (G, \cL', \cA')$ is {\em embedded} in $\xi$, if $L_i' \subset L_i$
for each $i$ and every operator $A_{ji}' = A_{ji}|_{L_i'}$ maps $L_i'$ to $L_j'$, whenever
$l_{ji} \in G$.
The embedding is strict if $L_i'$ is a proper subspace of $L_i$ for at least one~$i$. Thus,
an embedded triplet has the same multigraph and smaller spaces at the vertices.
A triplet $\xi = (G, \cL, \cA)$ is {\em reducible} if it has a strictly embedded triplet. Otherwise, it is called irreducible.
\smallskip

\noindent \textbf{Theorem~A}~\cite{CGP}. {\em An irreducible discrete system on graph with $\sigma =0$
possesses an invariant multinorm}. 
\medskip

\bigskip

\begin{center}
\textbf{9.2 Proofs  of Theorems~\ref{th.20},~\ref{th.30}, and~\ref{th.40}}
\end{center}
\bigskip

We begin with the following auxiliary fact, whose proof is placed in Appendix:  
\begin{lemma}\label{l.20}
Let $A_j, B_j$ be linear operators in~$\re^d$ and all~$B_j$ be invertible, $j=1, \ldots, n$. 
Suppose, for some segment $[\alpha , \beta]$, all operators from the 
set~$\bigl\{\, e^{tA_j}B_j, \ t\in [\alpha, \beta], \ j = 1, \ldots , n\bigr\}$
share a proper common invariant subspace;  then all~$A_j$ and $B_j$ share the same subspace. 
\end{lemma}

{\tt Proof of Theorem~\ref{th.30}.} 
 We consider the multigraph with vertices~$g_i$
associated to spaces~$L_i$. The set of  edges/operators
from~$L_i$ to $L_j$ is $\cA_{ji} \, = \, \bigl\{e^{tA_j} B_j, \ t \in [0, M-m] \bigr\}$. 
Clearly this is a compact set. If $\cA$ is irreducible, 
then by Lemma~\ref{l.20}, 
the set~$\bigl\{\, e^{tA_j}B_j, \ t\in [0, M-m], \ j = 1, \ldots , n\bigr\}$ is also irreducible.  This means that the discrete system on graph~$\xi$
is irreducible. Hence, we can invoke Theorem~A and obtain an invariant 
multinorm~$f$ for~$\xi$. By the definition of invariant multinorm, 
  for every point~$\bx$ from the unit sphere~$\cS_i$, 
  all the points~$A_{ji}\bx, \, A_{ji} \in \cA_{ji}$ are inside the unit ball~$\cB_j$ 
  and at least  one of them is on the sphere~$\cS_i$. 
  Therefore,  for every point~$\bx \in \cS_i$, 
  all the points~$B_je^{tA_j}\bx$ are inside~$\cB_j$ 
  and at least  one of them is on~$\cS_i$. Thus, $f$ is an invariant multinorm
  for the system~$\cA$, which completes the proof. 

{\hfill $\Box$}
\medskip

 {\tt Proof of Theorem~\ref{th.20}.} If the system is stable, then~$\sigma<0$. 
 Take arbitrary~$\alpha \in (\sigma, 0)$ and  
 an arbitrary multinorm~$\|\cdot \|\, = \, \bigl( \|\cdot \|_1\, , \cdots ,  \, \|\cdot \|_n\bigr)$ on~$(L_1, \ldots , L_n)$. 
 For every~$\bx \in L_i$, we 
 set~$f(\bx) \, = \, \sup\limits_{\by(\cdot), t \ge 0} e^{-\alpha t}\|\by(t)\|$, 
 where the supremum is computed over all trajectories~$\by$ such that~$\by(0) = \bx$. 
 Since $\alpha > \sigma$, it follows that the supremum is finite and 
 so~$f(\bx) < \infty$. Obviously~$f(\bx) > 0$ and $f$ is positively homogeneous. 
 Finally, $\by(t) \, = \, \Pi(t)\bx$, where $\Pi(t) \in \cE$ is the matrix product corresponding to the  switching law of the trajectory~$\by(\cdot)$. Hence, $f(\bx)$ is the supremum of 
 functions $e^{-\alpha t}\|\Pi(t)\bx\|$ over all proper matrices~$\Pi(t)$
and over all positive~$t$. Consequently, $f$ is convex as a supremum of convex functions. 
Thus, $f$ is a norm which possesses the 
property: for every trajectory~$\by(t)$ starting at~$\bx$, we have~$f\bigl(\by(t) \bigr) \le 
e^{\alpha t} f\bigl(\bx \bigr)$. Now consider an arbitrary trajectory~$\by(t)$
with switching points~$\{t_k\}_{k\in \n}$. Taking~$\bx = \by(t_k)$ and $t = t_{k+1}$
we get~$f\bigl(\by(t_{k+1}) \bigr) \le 
e^{\alpha t_{k+1}} f\bigl(\by(t_k) \bigr)\, < \, f\bigl(\by(t_{k}) \bigr)$, which proves that 
$f$ is a Luapunov norm. 

{\hfill $\Box$}
\medskip

\bigskip

{\tt Proof of Theorem~\ref{th.40}.} After a suitable shift it can  be 
assumed that  $\sigma(S) = 0$. If~$\cA$ is irreducible, then 
Theorem~\ref{th.30} provides  an invariant norm~$f$. Hence, there
are infinite trajectories for which~$f\bigl(\bx(t_k)\bigr) = 1$ 
for all switching points~$t_k\, , k \in \n$. 
Infinitely many points $\bx (t_k)$ belong to one space $L_i$ and, 
due to compactness of the unit sphere,
there is a convergent  subsequence $\bx (t_{k_s})$ as $s \to \infty$.
 Take arbitrary positive $\varepsilon < 1$.
For every $\delta > 0$,  there exists $N = N(\delta)$ such that
$f\, \bigl(\bx (t_{k_{s+1}}) - \bx (t_{k_{s}})\bigr)  < \delta$, whenever $s > N$.
On the other hand, those points belong to one trajectory, 
hence $\bx (t_{k_{s+1}}) = \Pi\, \bx (t_{k_{s}})$
for some product~$\Pi = \Pi(t_{k_{s+1}} - t_{k_{s}}) \, \in \, \cE_0$. 
 Thus, $f\bigl(\, (\Pi - I)\bx(t_{k_s}) \bigr) < \delta$ 
 and therefore~$\, \bigl| \rho (\Pi) - 1\bigr|\, \le  \, \frac{\varepsilon}{2}$ 
 whenever $\delta$ is small enough. Since $\, \ln \alpha \, \ge \, 
 (\alpha - 1) \, -  \, \frac12 |\alpha - 1|^2 $ for all positive~$\alpha$, we have  
 $\sigma(\Pi)\, = \, \ln \rho(\Pi) \,  \ge  \, - \frac{\varepsilon}{2} \, - \,  \frac{\varepsilon^2}{8}\, 
 >  \, - \varepsilon$, which completes the proof in the irreducible case. 
 
If the system is reducible, then all its matrices 
admit a simultaneous block upper-triangular factorisation with irreducible blocks 
$\cA_1, \ldots , \cA_r$. Since the spectrum of a block upper-triangular matrix 
is a union of spectra of the blocks, it follows that \linebreak 
$\sigma (\cA) \, = \, \max\, \bigl\{\sigma (\cA_1), \ldots , \sigma (\cA_r)\bigr\}$, 
where $\cA_k$ is the Lyapunov exponent of the restricted system (with the same 
switching interval~$[m,M]$) corresponding to the~$k$th block. 
Hence, if Theorem~\ref{th.40} holds in each irreducible block, it holds 
for the whole system. 

  {\hfill $\Box$}
 \smallskip 
 
 \textbf{Extra parts 1)}:  Thus, for the time restricted systems, the asymptotics of 
any trajectory along the real time axis is equivalent to that 
along the sequence of switching points. 
Therefore, for the Lyapunov exponent, we have $\sigma(S)\, = \ 
\inf\, \bigl\{ \alpha \in \re \ \bigl| \ \|\bx(t_k)\|\,  \le \, Ce^{\, \alpha t_k} \, , \ k \in \n\bigr\}$, 
where $\{t_k\}_{k\in \n}$ are the switching points of the trajectory~$\bx(t)$ and  the infimum is computed over all trajectories of the system.  
\bigskip 

\begin{center}
\textbf{9.3. Proof of Theorem~\ref{th.50}}
\end{center}
\bigskip

  Let us denote by $\|\cdot \|_P$ the multinorm~$\bigl\{ \|\cdot \|_{P_j}\bigr\}_{j =1}^n$, 
  where the norm~$\|\cdot \|_{P_j}$ in~$L_j$ is generated by the polytope~$P_j$. 
  The same for the multinorm~$\|\cdot \|_{P^{(k)}_j}$ generated by the 
  polytope~$P_j^{(k)}$ after the~$k$th iteration.

\begin{prop}\label{p.75}
There exists~$k_0 \in \n$ such that for all~$k \ge k_0$, 
all the polytopes $P_j^{(k)} \, = \, {\rm co}_s\, 
\cV_j^{(k)}$ possess nonempty interiors. 
\end{prop}
{\tt Proof.} For a given~$j$, we have $P_j^{(k)} \subset P_j^{(k+1)}$ for all~$k$, 
therefore,  either $P_j^{(k)}$ possesses a nonempty interior for all 
sufficiently big~$k$, or they are all contained in some proper subspace~$L_j' \in L_j$. 
If this happens for some~$j$, then the family~$\cA_j$ is reducible, which is a contradiction.   

  {\hfill $\Box$}
\smallskip

\begin{lemma}\label{l.30}
Let~$\|\cdot\|$ be an arbitrary norm in~$\re^d$ and 
$\bx: \, [0,\tau]\to \re^d$ be a $C^2$-curve.  
Then, for every~$t \in [0,\tau]$, 
the distance from the point~$\bx(t)$ to the segment~$[\bx(0),\bx(\tau)]$
does not exceed~$\frac{\tau^2}{8}\, \|\ddot \bx\|_{C[0, \tau]}$. 
\end{lemma}
{\tt Proof.} It can be assumed that~$\bx(0) = 0$, denote also~$\bx(\tau) = \ba$. 
Let the distance from~$\bx(t)$ to the segment~$[0,\ba]$
be equal to~$r$ and be attained at the point~$\bx(\xi)$. 
Let $\by$ be the closest to~$\bx(\xi)$
point of that segment. Denote~$\bh = \bx(\xi) \, - \, \by$. Thus, $\|\bh\| = r$. 
The segment~$[0,\ba]$ does not intersect the 
interior of the ball of radius~$r$ centered at~$\bx(\xi)$. 
Therefore, by the convex separation theorem, there exists 
a linear functional~$\bp \in \re^d, \, \|\bp\|^*=1$, 
that vanishes on the segment~$[0, \ba]$ and 
$(\bp, \bh)\, = \, \|\bh\|$. Since $(\bp, \bh)\, = \, 
(\bp, \bx(\xi)) \, - \, (\bp, \by) \, = \, (\bp, \bx(\xi))$, we have 
$(\bp, \bx(\xi)) \, = \, r$.  
Define the 
function~$f(t) \, = \, \bigl( \bp, \bx(t)\bigr)$. 
The maximum of~$f(t)$  on the segment~$[0,\tau]$ is attained at~$t=\xi$, 
hence~$f'(\xi) = 0$. 
Without loss of generality we assume that~$\xi \le \frac12\, \tau $, otherwise interchange 
the ends of the segment~$[0,\tau]$. The Tailor expansion 
of~$f$ at the point~$\xi$ gives~$f(t)\, = \, f(\xi) \, + \, f'(\xi)(t-\xi)\, + \, 
+\frac12 \, f'(\eta) (t-\xi)^2$, where $\eta \in [t, \xi]$. 
Since~$f'(\xi) = 0$, we obtain for~$t = 0$: 
$f(0)\, = \, f(\xi)\, + \, \frac12 \, f''(\eta) \xi^2$. Thus, 
$$
f(\xi) \ - \ f(0) \ = \ \bigl( \bp\, , \, \bx(\xi) - \bx(0)\bigr)  \ = \ 
  \bigl( \bp, \bx(\xi) \bigr)\, = \, r.
 $$ 
 Therefore, 
  $$
   r \ \le \ \frac12 \, \bigl|f''(\eta)\bigr| \xi^2\, \le \, \frac{\tau^2}{8}\, 
   \bigl| \bigl(\bp\, , \, \ddot \bx(\eta)\bigr) \bigr|  \ \le \ \frac{\tau^2}{8}\,  \bigl\|
 \ddot \bx  \bigr\|\, ,    
$$
which completes the proof. 

  {\hfill $\Box$}
\smallskip 

\begin{cor}\label{c.30}
Under the assumptions of~Lemma~\ref{l.30}, 
if $\bx$ is a solution of the differential
equation $\dot \bx \, = \, A\bx$, then, for every~$t \in [0,\tau]$, 
the distance from the point~$\bx(t)$ to the segment~$[\bx(0),\bx(\tau)]$
does not exceed
 \begin{equation}\label{eq.est3}
 \frac{\frac{\tau^2}{8}\, \|A^2\| }{1 \, - \, \frac{\tau^2}{8}\, \|A^2\| }\, \max\, \bigl\{\|\bx(0)\|\, , \|\bx(\tau)\|\,   \bigr\}
 \end{equation} 
\end{cor}
{\tt Proof.} 
Since  $\ddot \bx = A\dot \bx \, = \, A^2\bx$, it follows that 
 $|f''(\eta)|\, = \,  \bigl| \bigl(\bp, A^2 \bx(\eta)\bigr)\bigr|\, \le \, 
\bigl\|A^2\bx(\eta) \bigr\|\, \le \, \bigl\|A^2\bigr\|\, \bigl\|\bx (\eta) \bigr\|$. 
Thus, 
 \begin{equation}\label{eq.est4}
r \ \le \  \frac{\tau^2}{8}\, \|A^2\| \, \|\bx (\eta)\| \, . 
 \end{equation}
Let $\by$ be the closest to~$\bx(\eta)$ point of the segment~$[\bx(0), \bx(\tau)]$. 
We have 
$$
\bigl\|\bx (\eta)\bigr\| \ \le\   \|\by\|\, + \, r \ \le \ 
\|\by\| \ + \ \frac{\tau^2}{8}\, \|A^2\| \, \|\bx (\eta)\|, 
$$
therefore, 
$$
\bigl\|\bx (\eta)\bigr\|\ \le \ \frac{\|\by\|}{1 \, - \, \frac{\tau^2}{8}\, \|A^2\| }\ = \ 
\frac{\max\, \bigl\{\|\bx(0)\|\, , \|\bx(\tau)\|\,   \bigr\}}{1 \, - \, \frac{\tau^2}{8}\, \|A^2\| }
$$
The latter inequality follows from the convexity of the norm. 
Combining with~(\ref{eq.est4}), we complete the proof. 

  {\hfill $\Box$}
\smallskip

\begin{prop}\label{p.80}
If $\|\cdot \|$ is an arbitrary norm in~$\re^d$ and  $\bx(t)$ and is a solution of the differential
equation $\dot \bx \, = \, A\bx$, then, for every~$t \in [0,\tau]$,
we have 
 \begin{equation}\label{eq.est3.5}
 \bigl\|\bx(t) \bigr\|\quad \le \quad 
 \frac{1}{1 \, - \, \frac{\tau^2}{8}\, \|A^2\| }\, \max\, \Bigl\{\, \|\bx(0)\|\, , \, \|\bx(\tau)\|\,   \Bigr\}
 \end{equation} 
\end{prop}
{\tt Proof.} If $\by$ is the closest to~$\bx(t)$ point of 
the segment~$[\bx(0), \bx(\tau)]$, then 
$\|\bx(t)\|\le \|\by\| \, + \, \|\by - \bx(t)\|$. 
The convexity of the norm implies that~$\|\by\| \, \le \,   
\max\, \{\|\bx(0)\|, \|\bx(\tau)\|\}$. Estimating $\|\by - \bx(t)\|$
from above by~(\ref{eq.est3}) we arrive at~(\ref{eq.est3.5}).

  {\hfill $\Box$}
\smallskip

\begin{lemma}\label{l.40}\cite{P-bar}
There is a continuous function~$\psi(\delta, z)$ on~$\re_+^2$ such that 
$\psi(0, z) = 0$ for all~$z$ and for every $d\times d$  matrix~$A$, 
the following is true: if there is a vector $\bx$ 
such that $\|A\,\bx - \bx\| \, \le \, \delta \, \|\bx\|$, 
then~$A$ has an eigenvalue~$\lambda \in \co$ such that $|\lambda - 1| \, \le \, 
\psi (\delta, \|A\|)$.   
\end{lemma}

{\tt Proof of Theorem~\ref{th.50}.} {\em Algorithm terminates in Case~1}. 
We have $\mu > \delta$. Since the set of products~$\Pi(t_i, t)$
is contained in~$\cE_0$ and $\mu$ is the maximum of the values 
$\frac{1}{t-t_i}\sigma\bigl(\Pi(t_i, t)\bigr)$ along this set, is follows from 
Proposition~\ref{p.50} that~$\sigma(S)\ge \mu > - \delta$.
\smallskip

{\em  Algorithm does not terminates in Case~1}. Let us show that the algorithm terminates 
in Case~2, $\sigma(S) < \nu$, and $\|\cdot\|_P$ is the corresponding 
Lyapunov norm. 
\smallskip 

\noindent {\em Termination}. 
If, to the contrary,  Algorithm~1 does not halt in Case~2, 
then there is an infinite  trajectory on
the graph, for which all the switching points are not absorbed by the 
corresponding polytopes. This means that the switching point~$\bx(t_k)$, 
which appears in the~$k$th iteration on some~$L_j$, does not belong to the 
polytope~$P_j^{(k)}$. By Proposition~\ref{p.75}, 
there exists~$k_0 \in \n$ such that all $P_j^{(k_0)}$ possess nonempty interiors. 
Let they all contain a ball of radius~$r$ centered at zero.  
Since all those polytopes can only enlarge as $k$ grows, it follows that 
for all~$k\ge k_0$, the polytopes~$P_j^{(k)}$ contain this ball. 
Hence, $\|\bx(t_k)\| \ge r$ for all~$k \ge k_0$. Note that this trajectory 
is generated by the discrete system~$\cA_{\tau}$ on the graph. 
Hence, this system possesses a trajectory that does not converge to zero. 
Therefore, its joint spectral radius does not exceed one~\cite{CGP}. 
By continuity, there exists~$\alpha\ge 0$
such that the joint spectral radius of the discrete system~$(\cA - \alpha I)_{\tau}$
is equal to one
(the system~$(\cA - \alpha I)_{\tau}$ obtained from the system~$\cA - \alpha I$ by the same discretization as $\cA_{\tau}$ is 
obtained from~$\cA$). Clearly, this system is also irreducible. 
Hence, all its products are uniformly bounded by some constant~$C$~\cite{CGP}. 
Consider the discrete invariant polytope algorithm applied 
for the system~$(\cA - \alpha I)_{\tau}$ with the same starting point~$\bx (0)$
as in Algorithm~1. Since the joint spectral radius is equal to one, 
it follows that the discrete invariant polytope algorithm does not converge within finite time. Hence, there exists a discrete trajectory~$\tilde \bx(t)$ with $\tilde \bx(0) = \bx(0)$
and with some switching points~$\tilde \bx(\tilde t_k)$ such that 
every switching point does not belong to the corresponding polytope~$\tilde P^{(k)}_j$. 
Now we use the main argument.  Consider the trajectory of the original (continuous) system~$\cA$
with the same switching points~$\tilde t_k$. Denote it by~$\by(t)$ and by
$Q_j^{(k)}$ the corresponding polytopes generated by Algorithm~1. 
We have~$\by(t) = e^{\alpha t}\tilde x(t)$ for all~$t \ge 0$, every 
vertex of~$Q^{(k)}_j$ is obtained from the corresponding vertex of~$\tilde P^{(k)}_j$ 
by~$e^{\alpha t}$, where~$t$ is the time when that vertex is generated. 
Since $\tilde x(\tilde t_k)$ does not belong to the interior of~$\tilde P^{(k)}_j$ and  the point~$\by(\tilde t_k)$  appears later than all vertices of the 
polytope~$Q^{(k)}_j$ (later at least by the discretization step~$\tau$), 
it follows that 
 $\by(\tilde t_k)$ does not belong to the interior of~$e^{\tau \alpha} Q^{(k)}_j$
 and hence does not belong to~$Q^{(k)}_j$. 
 Thus, Algorithm~1 produces an infinite trajectory~$\by (t)$. 
 Hence, at every iteration $k$, the value~$\mu_k$ is bigger than 
 or equal to the maximum of the values~$\frac{1}{t-\tilde t_i}\sigma\bigl(\Pi(\tilde t_i, t)\bigr)$ over all products~$\Pi(\tilde t_i, t) \in \cE_0$. 
 This is bigger than or equal to the maximum of the corresponding 
 values~$\frac{1}{t-\tilde t_i}\sigma\bigl(\tilde \Pi(\tilde t_i, t)\bigr)$
 for the discrete system. Let as recall that for the discrete system, 
 the norms of all products~$\tilde \Pi(\tilde t_i, t)$ are bounded above by some constant~$C$. 
 Hence~$\|\tilde x(t)\| \le C \|x(0)\|$ for all~$t$, and 
 therefore, the trajectory~$\tilde x(t)$ possesses some limit point~$\bz \in L_i$. 
 For this limit point, there are products~$\tilde \Pi_r$ such that~$\|\tilde \Pi_r \bz \, - \, 
 \bz\|\to 0 $ as~$r \to \infty$. Since all those products are uniformly bounded, 
 it follows from Lemma~\ref{l.40}, that the distance from the point~$\lambda = 1$ to 
 the closest eigenvalue of~$\tilde \Pi_r$ tends to zero, and hence 
 $\limsup_{r\to \infty}\sigma (\tilde \Pi_r)\, \ge \, 0$. 
 On the other hand, the corresponding product of~$\Pi_r$ of the system 
 $\cA$ is equal to~$e^{\alpha \tau_t \tilde \Pi_r}$, where $\tau_r$ is 
 the time interval of the product~$\Pi_r$. Consequently, $\limsup_{r\to \infty}\sigma ( \Pi_r)\, \ge \, 0$ and hence 
 $\mu \ge 0$.  So Algorithm~1 must stop in Case~1, which is a contradiction.

\noindent {\em The upper bound~$\sigma(S) < \nu $ and the Lyapunov norm property}. 

By the construction of the polytopes~$P_j$, for every~$j = 1, \ldots  , n$ and $q\ne j$, we have 
$e^{\, s \tau A_j}B_jP_q \, \subset \, {\rm \int }\, P_j, \ s = 0, \ldots , N$. 
Since the set of all admissible triples~$(j,q,s)$ is finite, 
if follows that there is $\varepsilon > 0$ such that 
$$
e^{\, s \tau A_j}B_jP_q \, \subset \, (1-\varepsilon)\, P_j, \qquad
 j = 1, \ldots , n; \ q \ne j\, , \ s = 0, \ldots , N.
$$
Hence~$\|e^{\, s \tau A_j}B_j\|_P \le 1-\varepsilon$. 
To simplify the notation we denote $\|\cdot\|_{P_j}\, = \, \|\cdot\|_{j}$. 
Thus, for an arbitrary point~$\bx \in P_q$, we have 
$$
\bigl\|e^{\, s \tau A_j}B_j\bx\bigr\|_j \ \le \ (1-\varepsilon)\, 
\bigl\|\bx\bigr\|_q \quad \mbox{ for all} \ s = 0, \ldots , N. 
$$
Now take and arbitrary~$t \in [0, M-m]$, say, $t \in [s\tau , (s+1)\tau]$
for some~$s \le N-1$. Invoking Proposition~\ref{p.80} for $A= A_j, \bx(0) = e^{\, s \tau A_j}B_j\bx$
and $\bx(\tau) = e^{\, (s+1) \tau A_j}B_j\bx$ and taking into account that 
both $\|\bx(0)\|_j$ and $\|\bx(\tau)\|_j$ do not exceed~$(1-\varepsilon)\,\|\bx\|_q$, 
we obtain $\max\, \bigl\{ \|\bx(0)\|_j, \|\bx(\tau)\|_j \bigr\}\, \le \, 
\, (1-\varepsilon)\|\bx\|_q$. Therefore 
 \begin{equation}\label{eq.est5}
\bigl\|e^{\, t \tau A_j}B_j\bx\bigr\|_j  \quad \le \quad  \frac{1-\varepsilon}{1 \, - \, \frac{\tau^2}{8}\, \|A_j^2\|_j }\,  \|\bx\|_q\, . 
 \end{equation} 
Replacing~$\|A_j^2\|_j$ by $\|\cA^2\| \, = \,  \max_{j} \|A_j^2\|_j$, we increase 
the right hand side of~(\ref{eq.est5}). This yields 
$$
\bigl\|e^{\, t \tau A_j}B_j\bx\bigr\|_j  \quad \le \quad  \frac{1-\varepsilon}{1 \, - \, \frac{\tau^2}{8}\, \|\cA^2\| }\,  \|\bx\|_q\, , \qquad 
t \in [0, M-m]\,  
$$
Consequently, for an arbitrary trajectory starting at some point~$\bx_0$
for its $k$th switching point~$t_k$, we have 
$$
\|\bx(t_k)\| \quad \le \quad  
(1-\varepsilon)^{k-1}\, \left(1 \, - \, \frac{\tau^2}{8}\, \|\cA^2\| \, \right)^{-(k-1)} \, \|\bx_0\|\, .
$$
On the other hand, $t_k \ge (k-1)m$. Computing the logarithm of the both parts 
and taking the limit as~$k \to \infty$ gives 
$$
\sigma(S) \ \le \ -\, (1-\varepsilon)\frac1m\, \ln\, \left(1\, - \, 
\frac{\tau^2\, \|\cA^2\|}{8} \right)\, . 
$$
Substituting~$\tau = \frac{M-m}{N}$ we arrive at the upper bound of~(\ref{eq.est1}).

{\hfill $\Box$}
\medskip

\begin{center}
\textbf{Appendix} 
\end{center}

{\tt Proof of Proposition~\ref{p.10}.} The inequality~$\sigma < 0$ implies that 
every trajectory tends to zero, and hence, the system is stable. We need to establish
the converse: the stability implies that~$\sigma < 0$.  By the compactness argument,  for every~$T>0$, 
the set~$\|\Pi(t)\|, \,  t \in [0,T]$,  is uniformly bounded above by some constant~$H(T)$ over 
all switching laws from~$\cE(T)$.  
Denote by~$U(T)$ the subset of the unit sphere in~$\re^d$ such that 
there is a trajectory starting on that set for which~$\|\bx(T)\| \ge 1$. 
Consider two possible cases: 
\smallskip 

1) for some~$T$, the set~$U(T)$ is empty. In this case 
there is~$q<1$ such that $\|\bx(T)\| < q$ for all trajectories~$\bx(\cdot)$ starting on the unit sphere. 
Then splitting every infinite trajectory~$\bx(t)$ to the time intervals of length~$T$, 
we obtain~$\|\bx(Tk)\| < q^{\,k}, \, k \in \n$. Consequently, $\|\bx(t)\| \, < \, H(T)q^{\,k}$ 
for all~$t \in [kT, (k+1)T]$ and so~$\sigma \, \le \, T^{-1}\ln q\, < \, 0$. 
\smallskip 

2) the sets $U(T)$ are non-empty for all~$T$. Since they are all  compact and 
form an embedded system, they have a non-empty intersection. 
For every point~$\bx_0$ from this intersection and for every~$T$, there exists 
a trajectory~$\bx(t)$ such that~$\bx(0) = \bx_0$ and  $\bx(T)\, \ge \, 1$. Invoking  the 
compactness of the set of all trajectories on a segment, we obtain a trajectory, for which 
$\bx(T)\, \ge \, 1$ for all~$T>0$, hence the system is unstable.

{\hfill $\Box$}
\medskip

{\tt Proof of Lemma~\ref{l.20}.} Let~$L$ by the common invariant subspace of 
all~$e^{tA_j}B_j, \ t\in [\alpha, \beta], \ j = 1, \ldots , n$, and 
let~$C_j = e^{\alpha A_j}B_i$.  We first show that the operators~$A_j, C_j, \ i = 1, \ldots , n$ share the subspace~$L$. This is obvious for~$C_j$
since they belong to our set of operators (for~$t=0$). 
If some~$A_k$ does not respect~$L$, then  there exists~$\bx\in L$
such that~$A_k\bx \notin L$. Let $h > 0$ be the distance from~$A_k\bx$
to~$L$ and set~$\by \in C_k^{-1}\bx$. Since~$L$ is an invariant subspace for~$C_k$, 
it follows that~$\by \in L$.  Furthermore, since~$e^{tA_k}C_k\by \in L$ for all~$t \in [0, \beta - \alpha]$, and $e^{tA_k} \, = \, I + tA_k + O(t^2)$ as $t\to 0$, 
the distance from the point $(I + tA_k)C_k\by$ to~$L$ is also $O(t^2)$. 
On the other hand, $(I + tA_k)C_k\by \, =  \, (I + tA_k)\bx \, = \, \bx + tA_k \bx$, 
hence,  this distance is equal to~$th$. Thus, $th \, = \, O(t^2)$
as $t\to 0$, which is a contradiction. Consequently,  $L$ is invariant with respect to 
$A_k$ and $C_k$. Therefore, it is invariant for $e^{-\alpha A_k}$ and 
hence, for~$e^{-\alpha A_k}C_k = B_k$, which concludes the proof.   

  {\hfill $\Box$}


\medskip 

 
 
 
 \begin{thebibliography}{NN}
 
 \bibitem{AJ}
N.~Athanasopoulos and R.~M.~Jungers,
\newblock {\em Combinatorial methods for invariance and safety of hybrid systems},
\newblock Automatica J. IFAC, 98 (2018), 130--140.


\bibitem{B}
N.\,E.~Barabanov,
\newblock {\em Lyapunov indicator for discrete inclusions, I--III},
\newblock Autom. Remote Control, 49 (1988), No 2, 152--157.


\bibitem{Basso}
C.~Basso.
\newblock {\em Switch-mode power supplies spice simulations and practical
  designs}.
\newblock McGraw-Hill, Inc., New York, NY, USA, 1 edition, 2008.


\bibitem{BW} 
M.~A. Berger and Y.~Wang,
\newblock{\em Bounded semigroups of matrices},
\newblock Linear Alg. Appl.,
166 (1992) 21-27.
\smallskip


\bibitem{BCM}
F.~Blanchini, D,~Casagrande and S.~Miani,
\newblock {\em Modal and transition dwell time computation in switching systems: a set-theoretic approach},
\newblock Automatica J. IFAC, 46(9):1477--1482, 2010.

\bibitem{BM}
F.\,Blanchini, S.\,Miani,
\newblock {\em A new class of universal Lyapunov functions for the control of
uncertain linear systems},
\newblock IEEE Trans. Automat. Control,   44  (1999), no 3, 641--647.
\smallskip


\bibitem{BS}
C.~Briat and A.~Seuret.
\newblock Affine characterizations of minimal and mode-dependent dwell-times
  for uncertain linear switched systems.
\newblock {\em IEEE Trans. Automat. Control}, 58(5):1304--1310, 2013.

\bibitem{Chesi1}
G.~Chesi and P.~Colaneri.
\newblock Homogeneous rational {L}yapunov functions for performance analysis of
  switched systems with arbitrary switching and dwell time constraints.
\newblock {\em IEEE Trans. Automat. Control}, 62(10):5124--5137, 2017.

\bibitem{chesi0}
G.~Chesi, P.~Colaneri, J.~C. Geromel, R.~Middleton, and R.~Shorten.
\newblock A nonconservative {LMI} condition for stability of switched systems
  with guaranteed dwell time.
\newblock {\em IEEE Trans. Automat. Control}, 57(5):1297--1302, 2012.



\bibitem{CGPS21} 
Y.~Chitour, N.~Guglielmi, V.Yu.~Protasov, M.~Sigalotti, 
\newblock {\em Switching systems with dwell time: computing the maximal Lyapunov exponent}, 
\newblock  Nonlinear Anal. Hybrid Syst. 40 (2021), Paper No. 101021, 21 pp



\bibitem{CMS2}
Y.~Chitour, P.~Mason, and M.~Sigalotti.
\newblock A characterization of switched linear control systems with finite
  {$L_2$}-gain.
\newblock {\em IEEE Trans. Automat. Control}, 62(4):1825--1837, 2017.

\bibitem{CGP}
A.~Cicone, N.~Guglielmi, and V.~Y. Protasov.
\newblock Linear switched dynamical systems on graphs.
\newblock {\em Nonlinear Anal. Hybrid Syst.}, 29:165--186, 2018.


\bibitem{DO}
M.~Dehghan and C-J.~Ong,
\newblock {\em Characterization and computation of disturbance invariant sets
             for constrained switched linear systems with dwell time
             restriction},
\newblock Automatica J. IFAC, 48(9):2175--2181, 2012. 

\bibitem{FM}
L.\,Fainshil and M.\,Margaliot,
 \newblock {\em A maximum principle for the stability analysis of positive bilinear control systems with applications to positive linear switched systems},
 \newblock SIAM J. Control Optim.  50  (2012),  no. 4, 2193–2215.
\smallskip

\bibitem{GC}
J.~C. Geromel and P.~Colaneri.
\newblock Stability and stabilization of continuous-time switched linear
  systems.
\newblock {\em SIAM J. Control Optim.}, 45(5):1915--1930, 2006.

\bibitem{Teel00}
R.~Goebel, R.~G. Sanfelice, and A.~R. Teel.
\newblock {\em Hybrid dynamical systems}.
\newblock Modeling, stability, and robustness.
\newblock Princeton University Press, Princeton, NJ, 2012.


\bibitem{GLP17}
N.~Guglielmi, L.~Laglia, and V.~Protasov.
\newblock Polytope {L}yapunov functions for stable and for stabilizable {LSS}.
\newblock {\em Found. Comput. Math.}, 17(2):567--623, 2017.

\bibitem{GP13}
N.~Guglielmi and V.~Protasov.
\newblock Exact computation of joint spectral characteristics of linear
  operators.
\newblock {\em Found. Comput. Math.}, 13(1):37--97, 2013.

\bibitem{GP16}
N.~Guglielmi and V.~Protasov.
\newblock Invariant polytopes of linear operators with applications to
  regularity of wavelets and of subdivisions.
\newblock {\em SIAM J. Matrix Anal. Appl.}, 37(1):18--52, 2016.

\bibitem{HespanhaMorse}
J.~Hespanha and S.~Morse.
\newblock Stability of switched systems with average dwell-time.
\newblock In {\em Proceedings of the 38th IEEE Conference on Decision and
  Control}, 1999.

\bibitem{IngallsSontagWang}
B.~Ingalls, E.~Sontag, and Y.~Wang.
\newblock An infinite-time relaxation theorem for differential inclusions.
\newblock {\em Proceedings of the American Mathematical Society},
  131(2):487--499, 2003.
  

\bibitem{K}
V.~Kozyakin.
\newblock The {B}erger-{W}ang formula for the {M}arkovian joint spectral
  radius.
\newblock {\em Linear Algebra Appl.}, 448:315--328, 2014.

  \bibitem{K2010}
T.~Kr{\"o}ger.
\newblock {\em On-Line Trajectory Generation in Robotic Systems: Basic Concepts
  for Instantaneous Reactions to Unforeseen (Sensor) Events}.
\newblock Springer Berlin Heidelberg, Berlin, Heidelberg, 2010.


\bibitem{Liberzon}
D.~Liberzon.
\newblock {\em Switching in systems and control}.
\newblock Systems \& Control: Foundations \& Applications. Birkh\"{a}user
  Boston, Inc., Boston, MA, 2003.


\bibitem{LHM}
D.~Liberzon, J.P.~Hespanha, A.S.~Morse
\newblock {\em Stability of switched systems: a Lie-algebraic condition}, 
\newblock Syst.\& Contr. Letters, 37 (1999), no 3, 117--122. 

\bibitem{LiberzonMorse}
D.~Liberzon and A.~S. Morse.
\newblock Basic problems in stability and design of switched systems.
\newblock {\em IIEEE Control Systems Magazine}, 19:59--70, 1999.

\bibitem{Mej20}
T.~Mejstrik.
\newblock Improved invariant polytope algorithm and applications, 
 ACM Trans. Math. Softw., ACM Transactions on Mathematical Software 46 (2020), 
 no 3, p. 1--26,  https://doi.org/10.1145/3408891. 
\smallskip 



\bibitem{Morse}
A.~S. Morse.
\newblock Supervisory control of families of linear set-point controllers. {I}.
  {E}xact matching.
\newblock {\em IEEE Trans. Automat. Control}, 41(10):1413--1431, 1996.
\smallskip 

\bibitem{MP1}
A.~P. Molchanov and E.~S. Pyatnitski\u{\i}.
\newblock Lyapunov functions that define necessary and sufficient conditions
  for absolute stability of nonlinear nonstationary control systems. {III}.
\newblock {\em Avtomat. i Telemekh.}, (5):38--49, 1986.

\bibitem{MP2}
A.~P. Molchanov and Y.~S. Pyatnitskiy.
\newblock Criteria of asymptotic stability of differential and difference
  inclusions encountered in control theory.
\newblock {\em Systems Control Lett.}, 13(1):59--64, 1989.

\bibitem{PhJ1}
M.~Philippe, R.~Essick, G.~E. Dullerud, and R.~M. Jungers.
\newblock Stability of discrete-time switching systems with constrained
  switching sequences.
\newblock {\em Automatica J. IFAC}, 72:242--250, 2016.

\bibitem{PhJ3}
M.~Philippe, G.~Millerioux, and R.~M. Jungers.
\newblock Deciding the boundedness and dead-beat stability of constrained
  switching systems.
\newblock {\em Nonlinear Anal. Hybrid Syst.}, 23:287--299, 2017.

\bibitem{P-bar}
V.Yu.\,Protasov, 
\newblock The Barabanov norm is generically unique, simple, and easily computed.
\newblock {\em SIAM J. Contr. Optim.}, to appear in SIAM J. Cont. Opt.  (2022), 	arXiv:2109.12159. 

\bibitem{SFS}
M.~Souza, A.~Fioravanti, and R.~Shorten.
\newblock Dwell-time control of continuous-time switched linear systems.
\newblock In {\em Proceedings of the 54th IEEE Conference on Decision and
  Control}, pages 4661--4666, 2015.


\bibitem{VdSS}
A.~van~der Schaft and H.~Schumacher.
\newblock {\em An introduction to hybrid dynamical systems}, volume 251 of {\em
  Lecture Notes in Control and Information Sciences}.
\newblock Springer-Verlag London, Ltd., London, 2000.

\bibitem{WRDV}
Y.~Wang, N.~Roohi, G.~Dullerud, and M.~Viswanathan.
\newblock Stability of linear autonomous systems under regular switching
  sequences, 
\newblock In {\em Proceedings of the 54th IEEE Conference on Decision and
  Control}, pages 5445--5450, 2015.

\bibitem{Xiang2}
W.~Xiang.
\newblock On equivalence of two stability criteria for continuous-time switched
  systems with dwell time constraint,  
\newblock {\em Automatica J. IFAC}, 54:36--40, 2015.

\bibitem{Xiang1}
W.~Xiang.
\newblock Necessary and sufficient condition for stability of switched
  uncertain linear systems under dwell-time constraint, 
\newblock {\em IEEE Trans. Automat. Control}, 61(11):3619--3624, 2016.


 
 \end{thebibliography}
 \end{document}